\pdfoutput=1
\RequirePackage{ifpdf}
\ifpdf 
\documentclass[pdftex]{sigma}
\else
\documentclass{sigma}
\fi

\numberwithin{equation}{section}

\newtheorem{Theorem}{Theorem}[section]
\newtheorem*{Theorem*}{Theorem}

 { \theoremstyle{definition}

 }
\input{tcilatex}

\begin{document}
\allowdisplaybreaks

\newcommand{\arXivNumber}{2210.08771}

\renewcommand{\thefootnote}{}

\renewcommand{\PaperNumber}{087}

\FirstPageHeading

\ShortArticleName{ADE Bundles over Surfaces}

\ArticleName{ADE Bundles over Surfaces\footnote{This paper is a~contribution to the Special Issue on Enumerative and Gauge-Theoretic Invariants in honor of Lothar G\"ottsche on the occasion of his 60th birthday. The~full collection is available at \href{https://www.emis.de/journals/SIGMA/Gottsche.html}{https://www.emis.de/journals/SIGMA/Gottsche.html}}}

\Author{Yunxia CHEN~$^{\rm a}$ and Naichung Conan LEUNG~$^{\rm b}$}

\AuthorNameForHeading{Y.X.~Chen and N.C.~Leung}

\Address{$^{\rm a)}$~School of Mathematics, East China University of Science and Technology,\\
\hphantom{$^{\rm a)}$}~Meilong Road 130, Shanghai, China}
\EmailD{\href{mailto:yxchen76@ecust.edu.cn}{yxchen76@ecust.edu.cn}}

\Address{$^{\rm b)}$~The Institute of Mathematical Sciences and Department of Mathematics,\\
\hphantom{$^{\rm b)}$}~The Chinese University of Hong Kong, Shatin, N.T., Hong Kong}
\EmailD{\href{mailto:leung@math.cuhk.edu.hk}{leung@math.cuhk.edu.hk}}

\ArticleDates{Received March 07, 2022, in final form October 21, 2022; Published online November 05, 2022}

\Abstract{This is a review paper about ADE bundles over surfaces. Based on the deep connections between the geometry of surfaces and ADE Lie theory, we construct the corresponding ADE bundles over surfaces and study some related problems.}


\Classification{53C10; 14H60}

\renewcommand{\thefootnote}{\arabic{footnote}}
\setcounter{footnote}{0}

\section{Introduction}
Studies of curves on surfaces is one of the oldest problems in mathematics,
starting from the intersection theory of configurations of lines on a plane.
It is also an important tool in understanding the geometry of surfaces, see
for example \cite{BPV,Caporaso-Harris}. Enumeration of curves on surfaces is an old subject which dates back to the 19th century.
It has received much renewed interest and had become a very active research area, probably started from the ground breaking paper by Yau
and Zaslow \cite{Yau-Zaslow-96} which gave a beautiful formula for the
number of curves on a K3 surface in terms of a modular form. Namely the
number of rational curves on a K3 surface of a given degree is given by the
coefficients of the generating function%
\begin{equation*}
q/\Delta (q) =\prod\nolimits_{m=1}^{\infty }(1-q^{m})^{-24},
\end{equation*}%
where $\Delta $ is the well-known modular form of weight 12. Beauville \cite{Beauville-99} (see also~\cite{Fantechi-Gottsche-vanStraten-99}) gave a~mathematical approach to this formula by interpreting the BPS count in terms of compactified Jacobians of curves. For primitive
classes, the Yau--Zaslow formula was proven by Bryan and the second author
\cite{Bryan-Leung-Survey,Bryan-Leung-99,Bryan-Leung-YZ,Lee-Leung-YZ-K3-nonprim,Lee-Leung-g=1} and the full conjecture was
proven later by Klemm, Maulik, Pandharipande and Scheidegger \cite{KMPS-Pf-YZ-formula} (see also~\cite{Pandharipande-Thomas-16}) via mirror symmetry for a Calabi--Yau threefold with a K3 fibration.

In \cite{Gottsche-conj-98}, G\"{o}ttsche gave an intriguing generalization
of the Yau--Zaslow formula which applies to any surface $X$ and any genus $g$
as long as the curve class $C$ is sufficiently ample. Concretely the number
of genus $g$ curves in $\vert C\vert $ passing through an
appropriate number $r$ of points is given as the coefficient of $q^{C \cdot (C-K)/2}$ in the following power series in $q$%
\begin{equation*}
B_{1}^{K^{2}}B_{2}^{C\cdot K}\big(DG_{2}\big)^{r}\big(D^{2}G_{2}\big) /\big(\Delta\big(D^{2}G_{2}\big)\big)^{\chi(O_{X}) /2},
\end{equation*}%
where $D=q\frac{\rm d}{{\rm d}q}$, $G_{2}(q) =-\frac{1}{24}+\sum_{k>0}\big(\sum_{d|k}d\big) q^{k}$ is the Eisenstein series and $B_{1}(q) $ and $B_{2}(q)$ are universal power
series in $q$. The universality of this amazing conjecture has been
solved by Kool--Shende--Thomas \cite{Kool-Shende-Thomas-11} and Tzeng \cite{Tzeng-12}
independently in 2011 and 2012.

There is also a refined curve counting defined by Block and G\"{o}ttsche
which is related to tropical geometry \cite{Block-Gottsche-16} and real
algebraic geometry \cite{Gottsche-Shende-14} as well.

In this article, we study very different aspects of the geometry of curves
on surfaces, namely the intricate relationships between configurations of
low degree curves, for instance lines, and representation theory of simple
Lie algebras. The most famous example is probably the 27 lines on a cubic
surface discovered by Cayley and Salmon in 1849 and their relationships with
the exceptional Lie algebra $E_{6}$ (see, e.g., \cite{Dolgachev-2012,Manin-86}).

The organization of this paper is as follows. In Sections~\ref{sec2}, \ref{sec3}, \ref{sec4} and~\ref{sec5}, we introduce the famous examples that reflect the deep
connections between the geometry of surfaces and Lie theory. Specifically, we study ADE singularities in Section~\ref{sec2}, cubic surfaces in Section~\ref{sec3}, del Pezzo surfaces in Section~\ref{sec4}  and  ADE surfaces in Section~\ref{sec5}. Based on Section~\ref{sec5} (ADE surfaces), Sections~\ref{sec6} and~\ref{sec7} are related to F/string theory duality. In Sections~\ref{sec8}, \ref{sec9}, \ref{sec10} and~\ref{sec11}, we introduce some results related to deformation of ADE bundles over surfaces with ADE singularities, based on Section~\ref{sec2} (ADE singularities). We study Cox rings of ADE surfaces in Section~\ref{sec12}. And the final Section~\ref{sec13} is a summary.

\section{ADE singularities vs ADE Lie theory}\label{sec2}

In this section, we will recall ADE surface singularities, Lie algebras
of ADE types and their fundamental represenations.
Another intricate relationship between geometry of surfaces and Lie theory is the McKay correspondence \cite{Reid-2}. The simplest
type of surface singularity $p\in X$ is called a simple singularity
(also called a rational double point, canonical singularity, Du Val singularity or ADE
singularity) \cite{BPV}. Locally it is given by the quotient $\mathbb{C}^{2}/\Gamma $
of $\mathbb{C}^{2}$ by a finite subgroup $\Gamma \subset {\rm SL}(2,\mathbb{C})$. The exceptional curve $C=\bigcup_{i=1}^{r} C_{i}$ of its minimal
resolution $\tilde{X}\rightarrow X$ is a union of smooth rational curves $C_{i}$'s satisfying $C_{i}\cdot C_{i}=-2$, i.e., $(-2) $-curves,
and the configuration can be described by its dual graph which is a simply-laced Dynkin diagram, i.e., of ADE type.

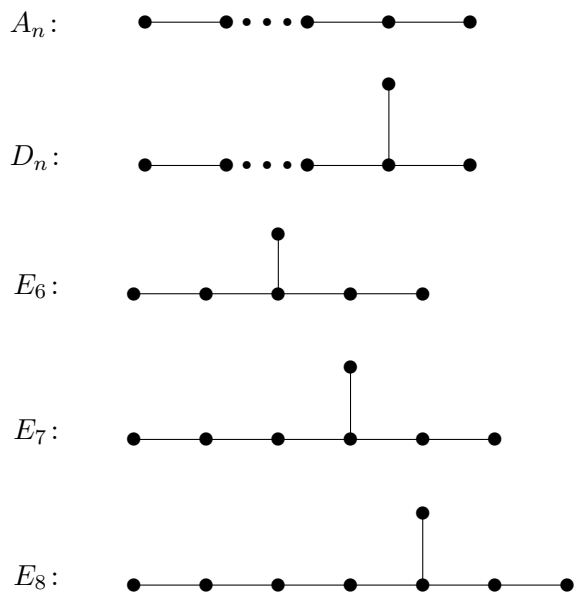
\begin{figure}[ht]\centering
\setlength{\unitlength}{0.9cm}\begin{picture}(5, 1.5) \put(-2,
0.85){$A_n\colon$}\put(0, 1){\circle*{.2}} \put(1.2, 1){\circle*{.2}} \put(2.4,
1){\circle*{.2}} \put(3.6, 1){\circle*{.2}} \put(4.8, 1){\circle*{.2}}
\put(0,1){\line(1, 0){1.2}} \put(1.5, 1){\circle*{.1}} \put(1.8,
1){\circle*{.1}}\put(2.1, 1){\circle*{.1}} \put(2.4, 1){\line(1, 0){1.2}}
\put(3.6, 1){\line(1, 0){1.2}} \end{picture}  \notag
\\[-5ex]
\setlength{\unitlength}{0.9cm}\begin{picture}(5, 3) \put(-2, 1){$D_n\colon$}
\put(0, 1){\circle*{.2}} \put(1.2, 1){\circle*{.2}}
\put(2.4,1){\circle*{.2}} \put(3.6,1){\circle*{.2}}
\put(4.8,1){\circle*{.2}} \put(3.6,2.2){\circle*{.2}} \put(1.5,
1){\circle*{.1}} \put(1.8, 1){\circle*{.1}} \put(2.1, 1){\circle*{.1}}
\put(0, 1){\line(1, 0){1.2}} \put(2.4, 1){\line(1, 0){1.2}} \put(3.6,
1){\line(1, 0){1.2}} \put(3.6, 1){\line(0, 1){1.2}} \end{picture}
\\[-5ex]
\setlength{\unitlength}{0.8cm}\begin{picture}(6, 3) \put(-2,
1){$E_6\colon$}\put(0, 1){\circle*{.2}} \put(1.2, 1){\circle*{.2}} \put(2.4,
1){\circle*{.2}} \put(3.6, 1){\circle*{.2}} \put(4.8, 1){\circle*{.2}}
\put(2.4,2.0){\circle*{.2}} \put(0, 1){\line(1, 0){1.2}} \put(1.2,
1){\line(1, 0){1.2}} \put(2.4, 1){\line(1, 0){1.2}} \put(3.6, 1){\line(1,
0){1.2}} \put(2.4, 1){\line(0, 1){1.0}} \end{picture}
\\[-3ex]
\setlength{\unitlength}{0.8cm}\begin{picture}(6, 3) \put(-2,
1){$E_7\colon$}\put(0, 1){\circle*{.2}} \put(1.2, 1){\circle*{.2}} \put(2.4,
1){\circle*{.2}} \put(3.6, 1){\circle*{.2}} \put(4.8, 1){\circle*{.2}}
\put(6, 1){\circle*{.2}} \put(3.6,2.2){\circle*{.2}} \put(0, 1){\line(1,
0){1.2}} \put(1.2, 1){\line(1, 0){1.2}} \put(2.4, 1){\line(1, 0){1.2}}
\put(3.6, 1){\line(1, 0){1.2}} \put(4.8,1){\line(1,0){1.2}} \put(3.6,
1){\line(0, 1){1.2}} \end{picture}
\\[-3ex]
\setlength{\unitlength}{0.8cm}\begin{picture}(6, 3) \put(-2,
1){$E_8\colon$}\put(0, 1){\circle*{.2}} \put(1.2, 1){\circle*{.2}} \put(2.4,
1){\circle*{.2}} \put(3.6, 1){\circle*{.2}} \put(4.8, 1){\circle*{.2}}
\put(6, 1){\circle*{.2}} \put(7.2,1){\circle*{.2}}
\put(4.8,2.2){\circle*{.2}} \put(0, 1){\line(1, 0){1.2}} \put(1.2,
1){\line(1, 0){1.2}} \put(2.4, 1){\line(1, 0){1.2}} \put(3.6, 1){\line(1,
0){1.2}} \put(4.8,1){\line(1,0){1.2}} \put(6,1){\line(1,0){1.2}} \put(4.8,
1){\line(0, 1){1.2}}
\end{picture}\vspace{-2ex}
\caption{Dynkin diagrams of ADE types.}\label{fig1}
\end{figure}

Recall that a simple Lie algebra\textsf{\ }$\mathfrak{g}$ is called
simply-laced if all roots have the same length and they are exactly those of
ADE types in the classification of simple Lie algebras \cite{Hum}, see Figure~\ref{fig1}.
Nodes in a Dynkin diagram label fundamental representations of $\mathfrak{g}$.

The Lie algebra $A_{n}=\mathfrak{sl}(n+1)$ is the algebra of
symmetries of a volume form on $V\simeq \mathbb{C}^{n+1}$. The fundamental
representations of $\mathfrak{sl}(n+1)$ consist of the standard
representation $V\simeq \mathbb{C}^{n+1}$ together with its exterior powers $\Lambda^{k}V$ with $k=2,3,\dots ,n$.

The Lie algebra $D_{n}=\mathfrak{o}(2n)$ is the algebra of infinitesimal
symmetries of a nondegenerate quadratic form $q\in S^{2}V^{\ast}$ on~$V\simeq \mathbb{C}^{2n}$. The fundamental representations of $\mathfrak{o}(2n)$
consist of the standard representation $V$ together with its exterior powers
$\Lambda^{k}V$ with $k=2,3,\dots,n-2$ and also the two spinor
representations $S^{+}$ and $S^{-}$.

\begin{figure}[ht]\centering
\setlength{\unitlength}{1cm}\begin{picture}(5, 1.5) \put(-1.55, 0.85){$A_n\colon$}\put(0,
1){\circle*{.2}} \put(1.2, 1){\circle*{.2}} \put(2.4, 1){\circle*{.2}}
\put(3.6, 1){\circle*{.2}} \put(4.8, 1){\circle*{.2}} \put(0,1){\line(1,
0){1.2}} \put(1.5, 1){\circle*{.1}} \put(1.8, 1){\circle*{.1}}\put(2.1,
1){\circle*{.1}} \put(2.4, 1){\line(1, 0){1.2}} \put(3.6, 1){\line(1,
0){1.2}}\put(-0.2,0.5){$V$}\put(0.8,0.5){${\wedge}^{2}V$}\put(1.8,0.5){${%
\wedge}^{n-2}V$}\put(3.2,0.5){${\wedge}^{n-1}V$}\put(4.7,0.5){${%
\wedge}^{n}V$} \end{picture}  \notag
\\[-2ex]
\setlength{\unitlength}{0.9cm}\begin{picture}(5, 3) \put(-2, 1){$D_n\colon$}
\put(0, 1){\circle*{.2}} \put(1.2, 1){\circle*{.2}}
\put(2.4,1){\circle*{.2}} \put(3.6,1){\circle*{.2}}
\put(4.8,1){\circle*{.2}} \put(3.6,2.2){\circle*{.2}} \put(1.5,
1){\circle*{.1}} \put(1.8, 1){\circle*{.1}} \put(2.1, 1){\circle*{.1}}
\put(0, 1){\line(1, 0){1.2}} \put(2.4, 1){\line(1, 0){1.2}} \put(3.6,
1){\line(1, 0){1.2}} \put(3.6, 1){\line(0, 1){1.2}}
\put(-0.2,0.5){$V$}\put(0.8,0.5){${\wedge}^{2}V$}\put(1.8,0.5){${%
\wedge}^{n-3}V$}\put(3.2,0.5){${\wedge}^{n-2}V$}\put(4.7,0.5){$S^{+}$}%
\put(3.8,2.1){$S^{-}$}\end{picture}
\\[-2ex]
\setlength{\unitlength}{0.8cm}\begin{picture}(6, 3) \put(-2,
1){$E_6\colon $}\put(0, 1){\circle*{.2}} \put(1.2, 1){\circle*{.2}} \put(2.4,
1){\circle*{.2}} \put(3.6, 1){\circle*{.2}} \put(4.8, 1){\circle*{.2}}
\put(2.4,2.0){\circle*{.2}} \put(0, 1){\line(1, 0){1.2}} \put(1.2,
1){\line(1, 0){1.2}} \put(2.4, 1){\line(1, 0){1.2}} \put(3.6, 1){\line(1,
0){1.2}} \put(2.4, 1){\line(0, 1){1.0}}
\put(-0.2,0.5){$V$}\put(4.7,0.5){$V^{*}$}\end{picture}
\\[-2ex]
\setlength{\unitlength}{0.8cm}\begin{picture}(6, 3) \put(-2,
1){$E_7\colon $}\put(0, 1){\circle*{.2}} \put(1.2, 1){\circle*{.2}} \put(2.4,
1){\circle*{.2}} \put(3.6, 1){\circle*{.2}} \put(4.8, 1){\circle*{.2}}
\put(6, 1){\circle*{.2}} \put(3.6,2.2){\circle*{.2}} \put(0, 1){\line(1,
0){1.2}} \put(1.2, 1){\line(1, 0){1.2}} \put(2.4, 1){\line(1, 0){1.2}}
\put(3.6, 1){\line(1, 0){1.2}} \put(4.8,1){\line(1,0){1.2}} \put(3.6,
1){\line(0, 1){1.2}} \put(-0.2,0.5){$V$}\end{picture}
\\[-2ex]
\setlength{\unitlength}{0.8cm}\begin{picture}(6, 3) \put(-2,
1){$E_8\colon $}\put(0, 1){\circle*{.2}} \put(1.2, 1){\circle*{.2}} \put(2.4,
1){\circle*{.2}} \put(3.6, 1){\circle*{.2}} \put(4.8, 1){\circle*{.2}}
\put(6, 1){\circle*{.2}} \put(7.2,1){\circle*{.2}}
\put(4.8,2.2){\circle*{.2}} \put(0, 1){\line(1, 0){1.2}} \put(1.2,
1){\line(1, 0){1.2}} \put(2.4, 1){\line(1, 0){1.2}} \put(3.6, 1){\line(1,
0){1.2}} \put(4.8,1){\line(1,0){1.2}} \put(6,1){\line(1,0){1.2}} \put(4.8,
1){\line(0, 1){1.2}} \put(-0.2,0.5){$V$} \end{picture}\vspace{-2ex}
\caption{Dynkin diagrams with fundamental representations labelled.}\label{fig2}
\end{figure}
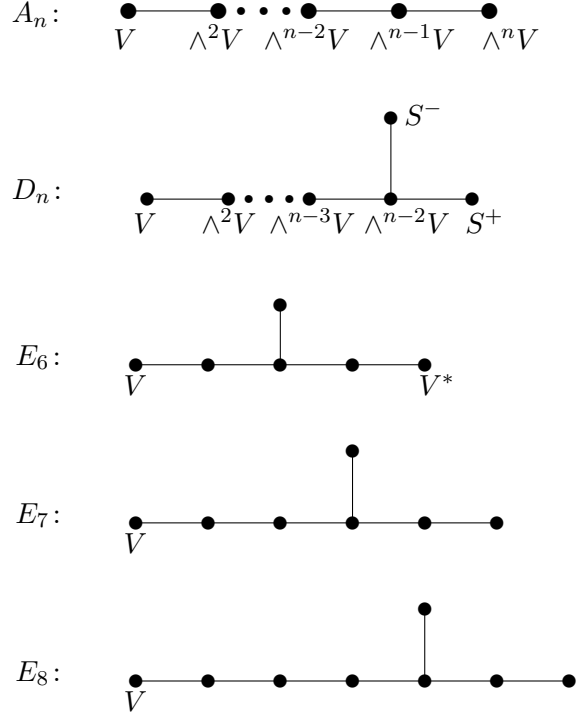

Furthermore, $E_{6}$ is the algebra of infinitesimal symmetries of a specific
cubic form $c\in S^{3}V^{\ast }$ on\ $V\simeq \mathbb{C}^{27}$ and $E_{7}$
is the algebra of infinitesimal symmetries of a specific quartic form $t\in
S^{4}V^{\ast }$ on\ $V\simeq \mathbb{C}^{56}$. The explicit forms of these
cubic form and quartic form can be described in terms of the geometry of del
Pezzo surfaces of degree 3 and 2 respectively. The situation for $E_{8}$ is
more complicated as its smallest representation $V$ is not a miniscule
representation, instead it is the adjoint representation of $E_{8}$. We will
call the above $V$'s as the standard representations~of~$\mathfrak{g}$, see Figure~\ref{fig2}.

\section[Cubic surfaces vs E\_6 Lie theory]
{Cubic surfaces vs $\boldsymbol{E_6}$ Lie theory}\label{sec3}

Let us come back to explain the relationships between cubic surfaces, or
more generally del Pezzo surfaces, with the representation theory of
exceptional Lie algebras of type $E$.

In this section, we will recall various ways to
realize the famous $27$ lines on cubic surfaces in $\mathbb{P}^{3}$.
Each such geometric setting corresponds to the branching rule for the $27$-dimensional standard representation of $E_{6}$ to a Lie subalgebra.
\begin{table}[h!]\centering\renewcommand{\arraystretch}{1.2}\small
\begin{tabular}{c|c}
\hline
Geometric settings & Lie subalgebra of $E_{6}$
\\ \hline
Degenerate to 3 planes & $\mathfrak{sl}(3) \times \mathfrak{sl}(3)\times \mathfrak{sl}(3)$
\\ \hline
Degenerate to plane + quadric surface & $\mathfrak{sl}(6) \times \mathfrak{sl}(2)$
\\ \hline
Blowing down a line & $\mathfrak{o}(10)$
\\ \hline
\end{tabular}
\end{table}

\noindent
We will also explain how to realize the $E_{6}$ structure from the
configuration of these 27 lines.

It is a classical result that there
are exactly 27 lines on a smooth cubic surfaces $X=\{f(z)\allowbreak=0\} \subset \mathbb{P}^{3}$ and the symmetry of their intersection
pattern is the Weyl group $W_{E_{6}}$ of $E_{6}$ (see, e.g., \cite{Manin-86,Reid}). Could we recover the Lie algebra $E_{6}$ itself, rather than
just its Weyl group? Indeed this can be done from the geometry of these 27
lines.

One way to locate these lines is to consider a family of cubic surfaces of
the form \cite{Segre}
\begin{equation*}
X(t) =\left\{ z_{1}z_{2}z_{3}+tf\left(
z_{0},z_{1},z_{2},z_{3}\right) =0\right\} \subset \mathbb{P}^{3}
\end{equation*}%
for a homogeneous polynomial $f$ of degree 3. Suppose $f$ is generic, then $%
X(t) $ is a smooth cubic surface for any generic $t\neq 0$. When
$t=0$, $X(0) $ is the union of 3 coordinate planes, which
certainly contain infinitely many lines. But only a few of them could
survive for nearby $t$. To~explain this, we note that the singular set of $%
X(0) $ is the union of 3 coordinate axes. In the smoothing $%
X\left( t\neq 0\right) $, each such singular point $p$ becomes a \textit{%
vanishing }$S^{1}$ as in the standard model of smoothing from $\left\{
xy=0\right\} $ to $\left\{ xy=t\right\} $, with the exception when $f$
vanishes at $p$. As~$\deg (f) =3$, $f$ vanishes at 3 points on
each coordinate axis. We call these \textit{unstable points} as without them
the total family of $X(t) $'s would be a semi-stable
degeneration for small $t$. A line on a coordinate plane, say $\left\{
z_{1}=0\right\} \subset X(0) $, survives on nearby $X(t) $'s if and only if it meets one of these 3 points on the $z_{2}$%
-axis, as well as on $z_{3}$-axis. Therefore the total number of lines on $%
X(0) $ which survive on nearby $X(t) $'s is $3\times
3+3\times 3+3\times 3=27.$

We will see that these 27 lines can be used to construct the fundamental
representation $L=\underline{27}$ of $E_{6}$. Here $\underline{27}$ means a
particular representation of $E_{6}$ of dimension 27. In fact, the above
description of these lines corresponds to the branching rule from $E_{6}$ to
its subalgebra $\mathfrak{su}(3)_{1}\times \mathfrak{su}(3)_{2}\times \mathfrak{su}(3)_{3}$:%
\begin{equation*}
\underline{27}=\underline{3_{1}}\times \underline{3_{2}^{\ast }}+\underline{%
3_{2}}\times \underline{3_{3}^{\ast }}+\underline{3_{3}}\times \underline{%
3_{1}^{\ast }}.
\end{equation*}%
Here $\underline{3_{i}}$ refers to the standard representation of the $i^{th}$-factor $\mathfrak{su}(3)_{i}$ in $\mathfrak{su}(3)_{1}\times \mathfrak{su}(3)_{2}\times \mathfrak{su}(3)_{3}$ and $\underline{3_{i}^{\ast }}$ is its dual representation.

Similarly, we could degenerate $X(t) $ to a union of a plane $H$
and a smooth quadric surface~$Q$ in $\mathbb{P}^{3}$, i.e., $X(0)
=H\cup Q$ \cite{LEUNG-Zhang-ADE-I}. Then there are $6$ points on the curve $C=H\cap Q$ which stay on~$X(t)$ for $t$ infinitesimally close to $0$. They play the same
role as the $3$ points on each coordinate axis in the previous example,
namely they are the unstable points for the family $X(t)$'s. A
line on the plane $H$ joining any $2$ of these $6$ points will deform to a
line on nearby $X(t)$'s. The total number of such lines is $\binom{6}{2}=15$. On the other hand, the quadric $Q$ has $2$ rulings, each is a~$\mathbb{P}^{1}$-family of lines. A line in $Q$ passing through one the these~$6$ points in $C$ will also deform to one in nearby $X(t) $, and
the total number of such lines is $6\times 2=12$. All together, we have $\binom{6}{2}+6\times 2=27$ lines on $X(t) $. This corresponds to
the branching rule for the fundamental representation $L_{E_{6}}=\underline{27}$ of $E_{6}$ to its subalgebra $\mathfrak{sl}(6) \times \mathfrak{sl}(2) \subset E_{6}$ as follows:%
\begin{equation*}
L_{E_{6}}\simeq \Lambda^{2}V_{6}+V_{6}\otimes V_{2},
\end{equation*}%
where $V_{6}\simeq \mathbb{C}^{6}$ and $V_{2}\simeq \mathbb{C}^{2}$ are the
standard representations of $\mathfrak{sl}(6) $ and $\mathfrak{sl}(2) $
respectively.

There is another way to see these lines once a particular line $l\subset X$
is given. Any hyperplane section containing $l$ must be of the form $l+C$
for some conic curve $C\subset X$ \cite{Manin-86,Reid}. The pencil of such hyperplane sections
degenerates into sum of three lines $l+l_{i}^{\prime }+l_{i}^{\prime \prime}$ five times $i=1,\dots,5$. These 10 divisors $l_{i}^{\prime }$'s and $l_{i}^{^{\prime \prime }}$'s determine the standard representation of $O(5,5) =D_{5}$ (see $D_{n}$-surfaces). The remaining 16 lines,
which do not intersect $l$, form a spinor representation of~$D_{5}$~and
\begin{equation*}
\underline{27}=\underline{1}+\underline{16}+\underline{10}
\end{equation*}%
is the branching rule from $E_{6}$ to $E_{5}=D_{5}$. From the surface
perspective, these three types of lines correspond to lines with
intersection numbers $-1$, $0$ and $1$ with the given line $l$.

We could realize the $E_{6}$ structure more concretely. We recall that the
standard representation~$\underline{27}$ of $E_{6}$ admits a cubic form%
\begin{equation*}
c\colon\ \underline{27}\otimes \underline{27}\otimes \underline{27}\rightarrow
\mathbb{C}
\end{equation*}%
such that $E_{6}\simeq \operatorname{Aut}(\underline{27},c)$, similar to
$D_{n}\simeq \operatorname{Aut}(\underline{2n},q)$ as the symmetry group of a
quadratic form on its standard representation. We consider the direct sum of
line bundles from lines on~$X$:%
\begin{equation*}
\mathcal{L}^{E_{6}}=\tbigoplus\limits_{i=1}^{27}O(l_{i}).
\end{equation*}%
If 3 of them $l_{1}$, $l_{2}$, $l_{3}$ form a \textit{triangle}, i.e., $l_{i}\cdot
l_{j}=1$ whenever $i\neq j$, then $l_{1}+l_{2}+l_{3}$ is a hyperplane
section. Therefore%
\begin{equation*}
O(l_{1}) \otimes O(l_{2}) \otimes O(l_{3}) \simeq O(1) \simeq K^{-1},
\end{equation*}%
where $K$ is the canonical line bundle of $X$. With suitable choices of
these isomorphisms \cite{LEUNG-Zhang-ADE-I} and using all triangles in $X$,
we obtain
\begin{equation*}
c_{\mathcal{L}}\colon\ \mathcal{L}^{E_{6}}\otimes \mathcal{L}^{E_{6}}\otimes
\mathcal{L}^{E_{6}}\rightarrow K^{-1}
\end{equation*}%
so that the fiberwise cubic structures give a conformal $E_{6}$-bundle $%
\mathcal{E}^{E_{6}}$ over the cubic surface $X$.

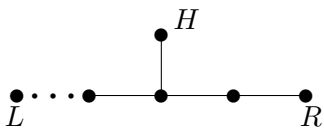
\begin{figure}[ht]\centering
\setlength{\unitlength}{0.8cm}\begin{picture}(6, 2.7)\put(0, 1){\circle*{.2}}
\put(1.2, 1){\circle*{.2}} \put(2.4, 1){\circle*{.2}} \put(3.6,
1){\circle*{.2}} \put(4.8, 1){\circle*{.2}} \put(2.4,2.0){\circle*{.2}}
\put(0.3, 1){\circle*{.1}} \put(0.6, 1){\circle*{.1}}\put(0.9,
1){\circle*{.1}} \put(1.2, 1){\line(1, 0){1.2}} \put(2.4, 1){\line(1,
0){1.2}} \put(3.6, 1){\line(1, 0){1.2}} \put(2.4, 1){\line(0, 1){1.0}}
\put(-0.2,0.5){$L$}\put(4.7,0.5){$R$}\put(2.6,2.1){$H$}\end{picture}\vspace{-2ex}
\caption{Dynkin diagram of $E_{n}$ with $L$, $R$ and $H$ labelled.}\label{fig3}
\end{figure}

The root lattice $\Lambda_{E_{6}}$ of $E_{6}$ can be identified as the
orthogonal complement $\langle K\rangle^{\bot }$ of $K$ in $\operatorname{Pic}(X) \simeq H^{2}(X,\mathbb{Z}) \simeq \mathbb{Z}^{1,6}$. Here $\mathbb{Z}^{1,6}$ denotes the lattice $\mathbb{Z}^{7}$ with
the quadratic form \mbox{$(1) \oplus (-1)^{\oplus 6}$},
which is isomorphic to $H^{2}(X,\mathbb{Z})$ equipped with the
intersection form. This is because any smooth cubic surface is a~blowup of $\mathbb{P}^{2}$ at 6 points. Note that $\alpha \in H^{2}(X,\mathbb{Z}) $ is a~root in~$\Lambda_{E_{6}}$ if and only if $\alpha \cdot
\alpha =-2$ and $\alpha \cdot K=0$, i.e., $\alpha \in \langle K\rangle^{\bot }\simeq \Lambda_{E_{6}}$. We denote the collection of all roots as $\Phi$. If $\alpha \in \Phi $ can be represented by an
effective divisor $C$, then $C\simeq \mathbb{P}^{1}$ is called an $(-2)$-curve. Furthermore, (1) $l\in \operatorname{Pic}(X) \simeq H^{2}(X,\mathbb{Z}) $ satisfying $l\cdot l=-1$ and $l\cdot K=-1$
is always represented by a unique line on $X$; (2) $f\in \operatorname{Pic}(X)
\simeq H^{2}(X,\mathbb{Z})$ satisfying $f\cdot f=0$ and $f\cdot
K=-2$ will be called a \textit{ruling, }or conic bundle, as it defines a $%
\mathbb{P}^{1}$-bundle $\Phi_{f}\colon X\rightarrow \vert f\vert
\simeq \mathbb{P}^{1}$ on $X$ with fiber degree 2; (3) $h\in \operatorname{Pic}(X) \simeq H^{2}(X,\mathbb{Z})$ satisfying $h\cdot h=1$
and $h\cdot K=-3$ gives $\Phi_{h}\colon X\rightarrow \vert h\vert
\simeq \mathbb{P}^{2}$ which realizes $X$ as a blowup of $\mathbb{P}^{2}$ at
6 points. Similar structures hold for all del Pezzo surfaces and they are
closely related to the fundamental representations~$L$,~$R$ and~$H$ of~$E_{n} $ corresponding the left node, the right node and the top node in the
Dynkin diagram of~$E_{n}$'s, see Figure~\ref{fig3}.

In fact, the Lie algebra $E_{6}$-bundle of fiberwise infinitesimal
automorphisms of $\big(\mathcal{L}^{E_{6}},c_{\mathcal{L}}\big)$ is
naturally%
\begin{equation*}
\mathcal{E}^{E_{6}}=O_{X}^{\oplus 6}\tbigoplus\limits_{\substack{ \alpha
\cdot K=0  \\ \alpha \cdot \alpha =-2}}O_{X}(\alpha),
\end{equation*}%
similar to the root space decomposition of Lie algebra $\operatorname{Lie}(E_{6}) =\mathfrak{h}\tbigoplus\limits_{\alpha \in \Phi }\mathfrak{g}_{\alpha}$.

Given any smooth anti-canonical curve $\Sigma \in \vert -K\vert$, then it is always an elliptic curve by the adjunction formula. For any $\alpha \in \Phi $, $O(\alpha) |_{\Sigma}$ has degree 0 and
therefore defines a flat line bundle over $\Sigma $. Indeed $\mathcal{E}^{E_{6}}|_{\Sigma }$ is a flat $E_{6}$-bundle over $\Sigma $ and every flat $E_{6}$-bundle over $\Sigma $ always arises this way for an embedding of $\Sigma $ into some cubic surface. Such a correspondence works for any $E_{n}$-bundle, or even any ADE bundle, as we will describe next.

\section[Del Pezzo surfaces vs E\_n Lie theory]
{Del Pezzo surfaces vs $\boldsymbol{E_n}$ Lie theory}\label{sec4}
Cubic surfaces are degree three del Pezzo surfaces. In this section, we
explain the analogous results between the configuration of lines in
degree $9-n$ del Pezzo surfaces and representation theory for the Lie
algebra $E_{n}$ for any $n$ between $1$ and $8$, as given in \cite{LEUNG-Zhang-ADE-I}.

A smooth surface $X$ is called a del Pezzo surface if $K^{-1}$ is ample. Its
degree $\deg X=K\cdot K$ is always between $1$ and $9$. A del Pezzo surface
of degree $d$ is always a blowup of $\mathbb{P}^{2}$ at $n=9-d$ points in
general position, with the exception of $\mathbb{P}^{1}\times \mathbb{P}^{1}$. We will denote them as $X_{n}$ and $X_{1}^{\prime }=\mathbb{P}^{1}\times
\mathbb{P}^{1}$ respectively. For example $X_{0}=\mathbb{P}^{2}$, $X_{5}$ is
the complete intersection of two quadric hypersurfaces in $\mathbb{P}^{4}$, $X_{6}$ is a cubic surface in $\mathbb{P}^{3}$, $X_{7}$ is a branched cover of $\mathbb{P}^{2}$ branched along a~quartic curve and $X_{8}$ has a pencil of
elliptic curves, i.e., $\vert -K\vert \simeq \mathbb{P}^{1}$. If
we blowup the unique base point of this pencil, we obtain a surface $X_{9}$
which is a rational elliptic surface. We will see in later sections that the
geometry of this surface is related to affine Kac--Moody Lie algebras of type
$\hat{E}_{8}$'s.

In fact, the $X_{n}$'s and $X_{1}^{\prime }$ correspond to simple Lie algebras
of type $E_{n}$, $1\leq n\leq 9$ and $\mathbb{C}$ if the Abelian Lie algebra
$\mathbb{C}$ is also included as a rank one simple Lie algebra. We recall
that $E_{n}$'s with $n\leq 5$ coincide with certain classical Lie algebras,
as can be seen directly from Dynkin diagrams. Concretely, we have $%
E_{5}\simeq D_{5}=\mathfrak{o}(10) $, $E_{4}\simeq A_{4}=\mathfrak{sl}(5)$, $E_{3}\simeq A_{2}\times A_{1}=\mathfrak{sl}(3) \times \mathfrak{sl}(2)$, $E_{2}\simeq \mathfrak{sl}(2) \times \mathbb{C}$, $E_{1}\simeq \mathbb{C}$, $E_{1}^{\prime }\simeq \mathfrak{sl}(2) $ and $E_{0}=0$. We remark that
there is a physical explanation of having two $E_{n}$'s when $n=1$ \cite{Julia} via supergravity in eleven dimensions where $E_{n}$ arises as one
compactifies the 11-dimensional space-time along a $n$-dimensional torus.

The first relationship between the geometry of a del Pezzo surface $X_{n}$
and the Lie algebra~$E_{n}$ is the following: the orthogonal complement $%
\langle K\rangle^{\bot}\subset H^{2}(X_{n},\mathbb{Z}) \simeq \mathbb{Z}^{n,1}$ is a root system $\Lambda_{E_{n}}$ for~$E_{n}$.

Similar to the cubic surface case, a line $l\subset X_{n}$ means a curve in $%
X_{n}$ of degree one with respect to the anti-canonical class, which is
equivalent to $l\in H^{2}(X_{n},\mathbb{Z})$ satisfying $l\cdot
l=-1$ and $l\cdot K=-1$. It is because each such cohomology class can be
represented as the Poincare dual of a unique line in $X_{n}$. The number of
lines in $X_{n}$ is equal to the dimension of the standard representation $L$
of $E_{n}$, unless $n=8$ in which case it is $\dim L-8=240$ because $L$ is
no longer a~minuscule representation. Explicitly these numbers are 1, 3, 6,
10, 16, 27, 56 and 240. Similarly $f\in \operatorname{Pic}(X) \simeq
H^{2}(X,\mathbb{Z})$ satisfying $f\cdot f=0$ and $f\cdot K=-2$
will be called a \textit{ruling}, or conic bundle, as it defines a $\mathbb{P}^{1}$-bundle $\Phi_{f}\colon X\rightarrow \vert f\vert \simeq \mathbb{P}^{1}$ on $X$ with fiber degree 2. The number of rulings in $X_{n}$ is equal to the dimension of the fundamental representation $R$ of $E_{n}$ for $n<7$.

We define%
\begin{gather*}
\mathcal{E}^{E_{n}} =O^{\oplus n}\tbigoplus\limits_{\substack{ \alpha
\cdot K=0  \\ \alpha \cdot \alpha =-2}}O(\alpha),
\qquad
\mathcal{L}^{E_{n}} =\tbigoplus\limits_{\substack{ l\cdot K=-1  \\ l\cdot
l=-1}}O(l),
\qquad
\mathcal{R}^{E_{n}} =\tbigoplus\limits_{\substack{ f\cdot K=-2  \\ f\cdot
f=0}}O(f).
\end{gather*}%
Then (1) $\mathcal{E}^{E_{n}}$ is always an $E_{n}$-bundle over $X_{n}$; (2)
$\mathcal{L}^{E_{n}}$ is an $\mathcal{E}^{E_{n}}$-representation bundle over~$X_{n}$ corresponding to the $E_{n}$ representation $L$, for $n<8$ and (3) $\mathcal{R}^{E_{n}}$ is an $\mathcal{E}^{E_{n}}$-representation bundle over $X_{n}$ corresponding to the $E_{n}$-representation $R$, for $n<7$. The restriction on $n$ is related to the fact that $L$ for $E_{8}$ and $R$ for $E_{7}$ are the adjoint representations.

If $l_{1}$ and $l_{2}$ are two lines in $X_{n}$ satisfying $l_{1}\cdot
l_{2}=1$, then $l_{1}+l_{2}$ is a ruling. This is reflected by the fact that
$R$ is an irreducible component of the tensor product $E_{n}$-representation
$L\otimes L$ and in fact we have a natural bundle homomorphism%
\begin{equation*}
\mathcal{L}^{E_{n}}\otimes \mathcal{L}^{E_{n}}\rightarrow \mathcal{R}^{E_{n}}
\end{equation*}%
over any $X_{n}$.

For example $X_{1}$ is the blowup of $\mathbb{P}^{2}$ and the exceptional
curve is its unique line. The two exceptional curves in $X_{2}$ together
with the strict transform of the line joining two blow up points in $\mathbb{P}^{2}$ give the 3 lines in $X_{2}$. Having two different types of lines in $X_{2}$ reflects that $L_{E_{2}}$ is a~reducible representation of $E_{2}$ and
it is also responsible for the fact that there are two different degree 8 del
Pezzo surfaces, namely $X_{1}$ and $X_{1}^{\prime }=\mathbb{P}^{1}\times
\mathbb{P}^{1}$. The $6=3\times 2$ lines in $X_{3}$ is reflecting the fact
that $X_{3}$ is the blow up of $\mathbb{P}^{2}$ in 2 different ways, which
in turns is the origin of the Cremona transformation. In terms of the
representation $L_{E_{3}}$ of $E_{3}\simeq \mathfrak{sl}(3) \times \mathfrak{sl}(2) $ we have%
\begin{equation*}
L_{E_{3}}\simeq V_{3}\otimes V_{2},
\end{equation*}%
where $V_{3}$ and $V_{2}$ are the standard representations of $\mathfrak{sl}(3) $ and $\mathfrak{sl}(2) $ respectively. On $X_{4}$, there are 10
lines and 5 rulings. Under the identification of $E_{4}$ with $\mathfrak{sl}(5) $, we have $R_{E_{4}}\simeq V_{5}$ the standard representation of $%
\mathfrak{sl}(5) $ and $L_{E_{4}}\simeq \Lambda^{3}V_{5}$. We could also
see the relationship $L_{E_{4}}\simeq \Lambda^{3}R_{E_{4}}$ from the fact
that every 3 distinct rulings on $X_{4}$ determines a unique line in $X_{4}$
which is a bisection to each of these 3 rulings. In fact, we have a natural
bundle isomorphism%
\begin{equation*}
\mathcal{L}^{E_{4}}\simeq \Lambda^{3}\mathcal{R}^{E_{4}}(K).
\end{equation*}

On $X_{5}$ there are 16 lines and 10 rulings, which are related to a spinor
representation $S_{10}^{+}$ and the standard representation $V_{10}$ of $%
E_{5}\simeq D_{5}=\mathfrak{o}(10) $. The defining quadratic form of $%
\mathfrak{o}(10) $ on $V_{10}\simeq \mathbb{C}^{10}$ is reflected by the
geometric fact that given any two rulings $f_{1}$ and $f_{2}$, we have $%
f_{1}\cdot f_{2}\leq 2$ and the equality sign holds if and only if $%
f_{1}+f_{2}=-K$. In fact, we have a fiberwise quadratic form $q$ over $X_{5}$,
\begin{equation*}
q\colon\ \mathcal{R}^{E_{5}}\otimes \mathcal{R}^{E_{5}}\rightarrow O(-K),
\end{equation*}%
so that the Lie algebra $E_{5}$-bundle over $X_{5}$ is the bundle of
infinitesimal symmetries of $q$ on~$\mathcal{R}^{E_{5}}$~\cite{Leung-ADE}. A~smooth $X_{5}$ is the complete intersection of two quadric hypersurfaces $Q_{0}$ and~$Q$ in~$\mathbb{P}^{4}$. If $Q_{0}$ varies in a one parameter
family $Q_{0}(t) $ and degenerates into a union of two
hyperplanes $Q_{0}(0) =H^{\prime }\cup H^{\prime \prime }$, thus
$X_{5}(0) =X_{5}^{\prime }\cup X_{5}^{\prime \prime }$, then its
geometry is governed by the reduction of the Lie algebra $E_{5}\simeq
D_{5}=\mathfrak{o}(10) $ to $A_{3}\times A_{1}\times A_{1}\simeq
D_{3}\times D_{2}=\mathfrak{o}(6) \times \mathfrak{o}(4)$. For instance
the branching rule%
\begin{equation*}
S_{10}^{+}\simeq S_{6}^{+}\otimes S_{4}^{+}+S_{6}^{-}\otimes S_{4}^{-}
\end{equation*}%
is reflected by the following geometric statements: $X_{5}^{\prime }$ and $X_{5}^{\prime \prime }$ are quadric surfaces and therefore each admits two
rulings. Lines on each ruling passing through one of the four points in the
curve $Q\cap H^{\prime }\cap H^{\prime \prime }$, which are the unstable
points for the total family $X_{5}(t) $, are exactly those lines
on $X_{5}(0) $ which will survive for nearby $X_{5}(t) $'s. Thus the $16=4\times 2+4\times 2$ lines on $X_{5}(t) $ is reflecting the above branching rule of $S_{10}^{+}$. The
relationship between the branching rule $V_{10}\simeq V_{6}+V_{4}$ and the
geometry of rulings on $X_{5}(t) $ can be described in a similar
fashion.

The relationship between the geometry of cubic surfaces $X_{6}$ and
representation theory of $E_{6}$ has been discussed in the previous section.
Recall that a smooth degree 2 del Pezzo surface $X_{7}$ is a double cover of
$\mathbb{P}^{2}$, branched along a quartic curve $Q\subset \mathbb{P}^{2}$.
A bitangent line to $Q$ in $\mathbb{P}^{2}$ determines a pair of lines in $X_{7}$ as its double cover, thus the 28 bitangents to $Q$ gives 56 lines in $X_{7}$. If $l$ and $l^{\prime }$ is any such pair of lines in $X_{7}$, then $l+l^{\prime }=-K$ and we have a fiberwise non-degenerate quadratic form $q$
on $\mathcal{L}^{E_{7}}$,
\begin{equation*}
q\colon\ \mathcal{L}^{E_{7}}\otimes \mathcal{L}^{E_{7}}\rightarrow O(-K).
\end{equation*}%
Furthermore, there is a fiberwise quartic form $f$ on $\mathcal{L}^{E_{7}}$,%
\begin{equation*}
f\colon\ \mathcal{L}^{E_{7}}\otimes \mathcal{L}^{E_{7}}\otimes \mathcal{L}%
^{E_{7}}\otimes \mathcal{L}^{E_{7}}\rightarrow O(-2K),
\end{equation*}%
so that its bundle of infinitesimal symmetries is the Lie algebra $E_{7}$%
-bundle $\mathcal{E}^{E_{7}}$ over $X_{7}$. We also have
\begin{equation*}
\mathcal{E}^{E_{7}}\simeq \bigl( O^{\oplus 7}+\mathcal{R}^{E_{7}}\bigr)(K).
\end{equation*}%
If we degenerate the smooth quartic branch curve $Q\subset \mathbb{P}^{2}$
to a double conic $2C$, then we obtain a~family of smooth del Pezzo surface $X_{7}(t) $ degenerating to a nonnormal surface $X_{7}(0) =X_{7}^{\prime}\cup X_{7}^{\prime \prime}$ with $X_{7}^{\prime}\simeq X_{7}^{\prime \prime }$ $\simeq \mathbb{P}^{2}\supset C$. There are
8 unstable points on $C\subset X_{7}(0) $. Each line in $X_{7}^{\prime}$ or~$X_{7}^{\prime \prime}$ passing through 2 of these 8
unstable points will survive on nearby $X_{7}(t) $'s, thus
giving $\binom{8}{2}+\binom{8}{2}=56$ lines on $X_{7}$. This is reflected by
the branching rule of $L_{E_{7}}$ of $E_{7}$ to its subalgebra $A_{7}=\mathfrak{sl}(8)$,
\begin{equation*}
L_{E_{7}}\simeq \Lambda^{2}V_{8}+\Lambda^{2}V_{8}^{\ast},
\end{equation*}%
where $V_{8}\simeq \mathbb{C}^{8}$ is the standard representation of $%
\mathfrak{sl}(8)$. For degree 1 del Pezzo surfaces $X_{8}$, the
relationships between their geometry with $E_{8}$ are discussed in \cite{Leung-ADE}.

\section{ADE surfaces vs ADE Lie theory}\label{sec5}

In this section, the above relationships between the representation theory of $E_{n}$ and the
geometry of del Pezzo surfaces will be generalized to other simply-laced Lie
algebras, i.e., Lie algebras of type ADE \cite{LEUNG-Zhang-ADE-I}. We call
the corresponding surfaces ADE surfaces. In this unified description, an $%
E_{n}$-surface is an $X_{n+1}$ together with a choice of a line, or
equivalently $X_{n}$ with a point in it.

First we describe $D_{n}$-surfaces. Given any ruling $f$ on $X_{n+1}$, its
linear system defines a~$\mathbb{P}^{1}$-bundle on $X_{n+1}$%
\begin{equation*}
\Phi_{f}\colon\ X_{n+1}\rightarrow \vert f\vert \simeq \mathbb{P}^{1}.
\end{equation*}%
For a generic $X_{n+1}$, there are $n$ singular fibers and each is a union
of two lines intersecting at a~point. To see this, we note the Euler
characteristic of $X_{n+1}$ equals
\begin{equation*}
\chi (X_{n+1}) =\chi \big(\mathbb{P}^{2}\big) +n+1=n+4.
\end{equation*}%
Smooth (resp. singular) fibers have Euler characteristic $2$ (resp.~$3$). We
have
\begin{equation*}
\chi (X_{n+1}) =\chi \big(\mathbb{P}^{1}\times \mathbb{P}^{1}\big) +\text{\#}\ (\text{singular fibers})
\end{equation*}%
and therefore there are $n$ singular fibers \cite{Bea}. Using these fiberwise lines,
i.e., $l\cdot f=0$, we construct a rank $2n$ vector bundle $\mathcal{L}^{D_{n}}$ over $X_{n+1}$,%
\begin{equation*}
\mathcal{L}^{D_{n}}=\tbigoplus\limits_{\substack{ l\cdot K=-1  \\ l\cdot
l=-1  \\ l\cdot f=0}}O(l).
\end{equation*}%
Notice that any two such lines $l$ and $l^{\prime}$ intersect if and only
if $l+l^{\prime}$ is a singular fiber of $\Phi_{f}$. When this happens, we have%
\begin{equation*}
O(l) \otimes O(l^{\prime}) \simeq O(f).
\end{equation*}%
Putting these isomorphisms together, we obtain a fiberwise quadratic form $q_{\mathcal{L}}$ on $\mathcal{L}^{D_{n}}$ over~$X_{n+1}$,
\begin{equation*}
q_{\mathcal{L}}\colon\ \mathcal{L}^{D_{n}}\otimes \mathcal{L}^{D_{n}}\rightarrow O(f).
\end{equation*}%
The bundle of infinitesimal symmetries of $\big(\mathcal{L}^{D_{n}},q_{\mathcal{L}}\big)$ is naturally%
\begin{equation*}
\mathcal{E}^{D_{n}}=O^{\oplus n}\tbigoplus\limits_{\substack{\alpha \cdot
K=0  \\ \alpha \cdot \alpha =-2  \\ \alpha \cdot f=0}}O(\alpha),
\end{equation*}%
which is a $D_{n}$-Lie algebra bundle over $X_{n+1}$ so that $\mathcal{L}^{D_{n}}$ is its representation bundle corresponding to the standard
representation of $D_{n}=\mathfrak{o}(2n)$. We call $(X_{n+1},f)$ a $D_{n}$-surface, denoted as~$X_{D_{n}}$.

In fact, the orthogonal complement of $K$ and $f$ in $H^{2}(X_{n+1},\mathbb{Z})$ is the root lattice of type $D_{n}$
\begin{equation*}
\langle K,f\rangle^{\bot }\simeq \Lambda_{D_{n}}.
\end{equation*}%
Furthermore, this holds true for any $n$ without the restriction $n\leq 8$.

In a similar fashion, an $A_{n}$-surface $X_{A_{n}}$ is a pair $(X_{n+1},h)$ with $h\cdot K=-3$ and $h\cdot h=1$. The rank $n+1$
vector bundle%
\begin{equation*}
\mathcal{L}^{A_{n}}=\tbigoplus\limits_{\substack{ l\cdot K=-1  \\ l\cdot
l=-1  \\ l\cdot h=0}}O(l)
\end{equation*}%
over $X_{A_{n}}$ admits a fiberwise determinant morphism%
\begin{equation*}
\det\colon\ \Lambda^{n+1}\mathcal{L}^{A_{n}}\overset{\sim}{\rightarrow}O(K+3h).
\end{equation*}%
This is because
\begin{equation*}
\Phi_{h}\colon\ X_{n+1}\rightarrow \vert h\vert \simeq \mathbb{P}^{2},
\end{equation*}%
which realizes $X_{n+1}$ as a blowup of $\mathbb{P}^{2}$ at $n+1$ points and
exceptional curves for $\Phi_{h}$ are precisely those lines used to
construct $\mathcal{L}^{A_{n}}$. The bundle of infinitesimal symmetries of $\big(\mathcal{L}^{A_{n}},\det\big)$ is naturally%
\begin{equation*}
\mathcal{E}^{A_{n}}=O^{\oplus n}\tbigoplus\limits_{\substack{ \alpha \cdot
K=0  \\ \alpha \cdot \alpha =-2  \\ \alpha \cdot h=0}}O(\alpha),
\end{equation*}%
which is an $A_{n}$-Lie algebra bundle over $X_{A_{n}}$ so that $\mathcal{L}^{A_{n}}$ is its representation bundle corresponding to the standard
representation of $A_{n}=\mathfrak{sl}(n+1)$. Again the orthogonal
complement of $K$ and $h$ in $H^{2}(X_{n+1},\mathbb{Z})$ is the
root lattice of type $A_{n}$%
\begin{equation*}
\langle K,h\rangle^{\bot }\simeq \Lambda_{A_{n}}.
\end{equation*}

We remark that an $E_{n}$-surface can be interpreted as $(X_{n+1},l)$, as in the $D_{n}$ and $A_{n}$ cases. But since $l$ is an
exceptional curve and therefore it can be blown down to $X_{n}$. We also
remark that a choice of $h$, $f$ or $l$ in $X_{n+1}$ defines an $A_{n}$-surface, $D_{n}$-surface or $E_{n}$-surface accordingly. From the Lie
theoretical perspective, they correspond to the fundamental representations $H$, $R$ and $L$ for $E_{n+1}$. By removing the corresponding nodes in the $E_{n+1}$ Dynkin diagram, we also obtain Dynkin diagrams of type $A$, $D$ and $E$ respectively.

\pagebreak
An ADE surface is a surface $X_{n+1}$, $\mathbb{P}^{2}$ blown up at $n+1$
points in general position, together with a divisor $h$, $r$ or $l$ as below:
\begin{table}[h!]\centering\renewcommand{\arraystretch}{1.2}\small
\begin{tabular}{c|c}
\hline
$A_{n}$\text{-surface} & $(X_{n+1},h) \text{ with }h\cdot K=-3\text{
and }h\cdot h=1$\\ \hline
$D_{n}$\text{-surface} & $(X_{n+1},r) \text{ with }r\cdot K=-2\text{
and }r\cdot r=0$  \\\hline
$E_{n}$\text{-surface} & $(X_{n+1},l) \text{ with }l\cdot K=-1\text{
and }l\cdot l=-1$  \\ \hline
\end{tabular}
\end{table}

\section{F/String theory duality}\label{sec6}

In this section,
we recall a physical motivation for the
construction of $E_{n}$-bundles over del Pezzo surfaces of degree $9-n$
given in \cite{Clingher-Morgan-05,Donagi-97,Donagi-Wijnholt-20,Friedman-Morgan-Witten}. Physically, if $G$ is a simple subgroup of $E_{8}\times E_{8}$, then $%
G $-bundles are related to the duality between F-theory and heterotic string theory.
Among other things, this duality predicts that the moduli of flat $E_{n}$-bundles over a fixed elliptic curve~$\Sigma $ can be identified with the
moduli of del Pezzo surfaces with a fixed anti-canonical curve $\Sigma$.

Given any smooth elliptic curve $\Sigma $ and any ADE Lie group $G$ of
adjoint type, if $\Sigma $ is an anti-canonical curve in a $G$-surface $X$,
then the natural $\mathfrak{g}$ Lie algebra bundle $\mathcal{E}^{\mathfrak{g}}$ restricts to a~flat $G$-bundle over $\Sigma $. In fact, every flat $G$-bundle on $\Sigma $ arises this way. To state the result, we denote $\mathcal{M}_{\Sigma}^{G}$ to be the moduli space of flat $G$-bundles over $\Sigma$ and $\mathcal{S}_{\Sigma}^{G}$ to be the moduli space of pairs $(X,\Sigma \in \vert -K_{X}\vert)$ with $X$ being an ADE surface of type $G$ and $\Sigma $ is an anti-canonical curve in $X$.

\begin{Theorem}[{\cite{LEUNG-Zhang-ADE-I}}]\label{th6.1}
Given $\Sigma $ and $G$ as above, there is an open dense embedding%
\begin{equation*}
\Phi\colon\ \mathcal{S}_{\Sigma }^{G}\rightarrow \mathcal{M}_{\Sigma }^{G}
\end{equation*}%
given by the restriction of the natural $\mathfrak{g}$-bundle $\mathcal{E}^{\mathfrak{g}}$ over $X$ to $\Sigma $. Furthermore, there is a natural
compactification $\mathcal{\bar{S}}_{\Sigma}^{G}$ of $\mathcal{S}_{\Sigma}^{G}$ given by those surfaces $X$ equipped with $G$-configurations and $\Phi$ extends to an isomorphism $\Phi\colon\mathcal{\bar{S}}_{\Sigma }^{G}\overset{\sim}{\rightarrow }\mathcal{M}_{\Sigma}^{G}$.
\end{Theorem}

This particular form is given in \cite{LEUNG-Zhang-ADE-I} and its
generalization for non-simply laced $G$ is given in~\cite{LEUNG-Zhang-ADE-II,Leung-Zhang-12}. Such a correspondence was originally motivated from the duality between
F-theory and heterotic string theory in physics by the work of
Friedman--Morgan--Witten~\cite{Friedman-Morgan-Witten} and Do\-nagi~\cite{Donagi-97} where different proofs of this correspondence are also given.

Let us very briefly describe how such a correspondence arises from the
physical duality between F-theory and heterotic string theory. The
space-time in F-theory is a Calabi--Yau fourfold~$Z$ equipped with an
elliptic K3 fibration over a complex surface $B$ and the space-time in
heterotic string theory is a Calabi--Yau threefold $Y$ equipped with an
elliptic fibration over the same complex surface $B$ and coupled with an $%
E_{8}\times E_{8}$ Hermitian Yang--Mills bundle over $Y$. When the two
theories are dual to each other, in a certain adiabatic limit, the duality
becomes a fiberwise duality. Namely an elliptic K3 fiber $X$ over $b\in B$
in $Z$ is dual to an elliptic curve fiber $\Sigma $ over $b\in B$ in $Y$
coupled with a flat $E_{8}\times E_{8}$-bundle over $\Sigma $. To obtain a
geometric correspondence, the dilaton field in the string theory should
vanish, which corresponds to a~type~II degeneration of the elliptic K3
surface, that is $X$ is a~fiber sum $X_{1}\#_{\Sigma }X_{2}$, each $X_{i}$
is a~rational elliptic surface with section with $\Sigma$ being a fiber. In
particular, $X_{i}$ is a blowup of $\mathbb{P}^{2}$ at~$9$ points with one
exceptional curve identified with the section, namely $X_{i}$ is an $E_{8}$-surface as defined previously. Hence each $(X_{i},\Sigma)$ is
an $E_{8}$-surface and therefore gives an $E_{8}$-bundle over~$\Sigma $ by
our above discussions. Together we obtain the flat $E_{8}\times E_{8}$-bundle over $\Sigma $.

\section{Generalization to Kac--Moody cases}\label{sec7}

In this section, we generalize the above results for $E_{n}$-bundles in Theorem~\ref{th6.1} to Kac--Moody $\hat{E}_{n}$-bundles in~\cite[Theorems~2 and~3]{Leung-Xu-Zhang-12}.
If we blowup $\mathbb{P}^{2}$ at 9 points, then the resulting surface $X_{9}$
is no longer a del Pezzo surface. As $K\cdot K=0$, $\langle
K_{X}\rangle^{\perp }\subset H^{2}(X_{9},\mathbb{Z})$ is
only non-positive definite and there are infinitely many roots. The latter
is because for any root $\alpha $, namely $\alpha \cdot \alpha =-2$ and $\alpha \cdot K=0$, $\alpha +nK$ is also a root for any integer $n$ \cite{Chen-Leung-16,Manin-86}.

In fact, $\Lambda_{E_{9}}\triangleq \langle K_{X}\rangle^{\perp}$ is the root system of the affine Kac--Moody Lie algebra $\hat{E}_{8}$, or
the loop algebra $LE_{8}$ of $E_{8}$, up to central extension. The following
is the list of Dynkin diagrams of simply-laced affine Kac--Moody Lie
algebras, which coincides with extended Dynkin diagrams of simply-laced
simple Lie algebras, namely ADE types (see \cite{Kac} for details):
\begin{figure}[ht]\centering
\setlength{\unitlength}{0.9cm}\begin{picture}(5, 2.5) \put(-2,
1){$\hat{A}_n\colon$}\put(0, 1){\circle*{.2}} \put(1.2, 1){\circle*{.2}}
\put(2.4, 1){\circle*{.2}} \put(3.6, 1){\circle*{.2}} \put(4.8,
1){\circle*{.2}} \put(2.4, 2.2){\circle*{.2}} \put(0,1){\line(2, 1){2.4}}
\put(4.8,1){\line(-2, 1){2.4}} \put(0,1){\line(1, 0){1.2}} \put(1.5,
1){\circle*{.1}} \put(1.8, 1){\circle*{.1}}\put(2.1, 1){\circle*{.1}}
\put(2.4, 1){\line(1, 0){1.2}} \put(3.6, 1){\line(1, 0){1.2}} \end{picture}
\notag
\\[-3ex]
\setlength{\unitlength}{0.9cm}\begin{picture}(5, 3) \put(-2,
1){$\hat{D}_n\colon$} \put(0, 1){\circle*{.2}} \put(1.2, 1){\circle*{.2}}
\put(2.4,1){\circle*{.2}} \put(3.6,1){\circle*{.2}}
\put(4.8,1){\circle*{.2}} \put(3.6,2.2){\circle*{.2}}
\put(1.2,2.2){\circle*{.2}} \put(1.2, 1){\line(0, 1){1.2}} \put(1.5,
1){\circle*{.1}} \put(1.8, 1){\circle*{.1}} \put(2.1, 1){\circle*{.1}}
\put(0, 1){\line(1, 0){1.2}} \put(2.4, 1){\line(1, 0){1.2}} \put(3.6,
1){\line(1, 0){1.2}} \put(3.6, 1){\line(0, 1){1.2}} \end{picture}
\\
\setlength{\unitlength}{0.8cm}\begin{picture}(6, 3) \put(-2,
1){$\hat{E}_6\colon$}\put(0, 1){\circle*{.2}} \put(1.2, 1){\circle*{.2}}
\put(2.4, 1){\circle*{.2}} \put(3.6, 1){\circle*{.2}} \put(4.8,
1){\circle*{.2}} \put(2.4,2.0){\circle*{.2}} \put(2.4,3.0){\circle*{.2}}
\put(0, 1){\line(1, 0){1.2}} \put(1.2, 1){\line(1, 0){1.2}} \put(2.4,
1){\line(1, 0){1.2}} \put(3.6, 1){\line(1, 0){1.2}} \put(2.4, 1){\line(0,
1){1.0}} \put(2.4,2.0){\line(0, 1){1.0}} \end{picture}
\\[-2ex]
\setlength{\unitlength}{0.8cm}\begin{picture}(6, 3) \put(-2,
1){$\hat{E}_7\colon$}\put(0, 1){\circle*{.2}} \put(1.2, 1){\circle*{.2}}
\put(2.4, 1){\circle*{.2}} \put(3.6, 1){\circle*{.2}} \put(4.8,
1){\circle*{.2}} \put(6, 1){\circle*{.2}} \put(7.2,1){\circle*{.2}}
\put(3.6,2.2){\circle*{.2}} \put(0, 1){\line(1, 0){1.2}} \put(1.2,
1){\line(1, 0){1.2}} \put(2.4, 1){\line(1, 0){1.2}} \put(3.6, 1){\line(1,
0){1.2}} \put(4.8,1){\line(1,0){1.2}} \put(6,1){\line(1,0){1.2}} \put(3.6,
1){\line(0, 1){1.2}} \end{picture}
\\[-1ex]
\setlength{\unitlength}{0.8cm}\begin{picture}(6, 3) \put(-2,
1){$\hat{E}_8\colon$}\put(0, 1){\circle*{.2}} \put(1.2, 1){\circle*{.2}}
\put(2.4, 1){\circle*{.2}} \put(3.6, 1){\circle*{.2}} \put(4.8,
1){\circle*{.2}} \put(6, 1){\circle*{.2}} \put(7.2,1){\circle*{.2}}
\put(8.4,1){\circle*{.2}} \put(6,2.2){\circle*{.2}} \put(0, 1){\line(1,
0){1.2}} \put(1.2, 1){\line(1, 0){1.2}} \put(2.4, 1){\line(1, 0){1.2}}
\put(3.6, 1){\line(1, 0){1.2}} \put(4.8,1){\line(1,0){1.2}}
\put(6,1){\line(1,0){1.2}} \put(7.2,1){\line(1,0){1.2}} \put(6, 1){\line(0,
1){1.2}} \end{picture}\vspace{-2ex}
\caption{Dynkin diagrams of affine ADE types.}\label{fig4}
\end{figure}
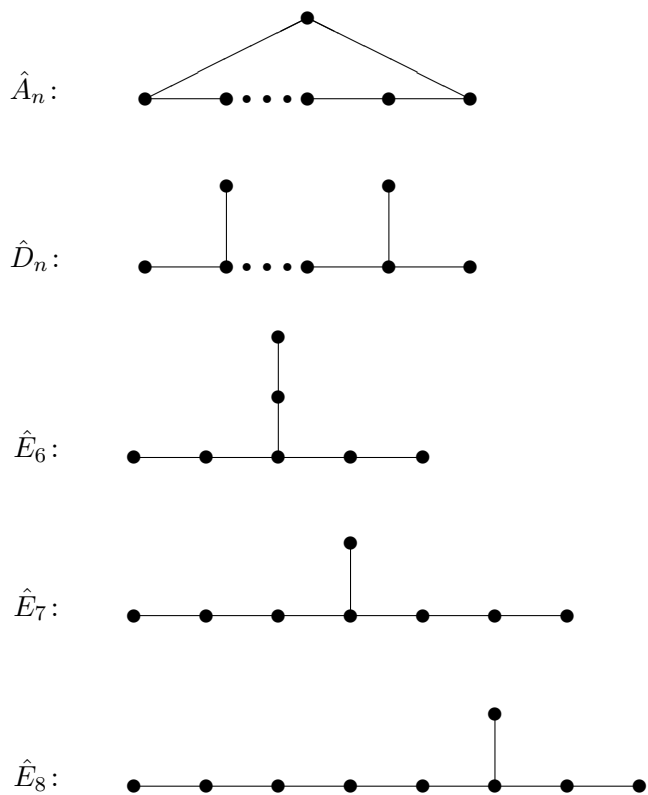

We could also realize other $\hat{E}_{n}$'s with $n\leq 8$ from rational
surfaces as follows: Let $C_{1},\dots ,C_{8-n}$ be an $A_{8-n}$-chain of $%
\left( -2\right) $-curves in $X_{9}$, i.e., $C_{i}\simeq \mathbb{P}^{1}$, $%
C_{i}\cdot C_{i}=-2$, $C_{i}\cdot C_{i+1}=1$ and other $C_{i}\cdot C_{j}$'s
are zero. Then we could blow down this chain of $\left( -2\right) $-curves
in $X_{9}$ to obtain a rational surface, denoted $X_{\hat{E}_{n}}$, with a
rational singular point of type $A_{8-n}$. Equivalently $X_{\hat{E}_{n}}$ is
a singular del Pezzo surface of degree $8-n$ with a canonical singularity of type $A_{8-n}$, for $n \leq 7$. When $n=8$,
further care is needed as the representation $L_n$ is no longer a miniscule representation in this case. We call these $X_{%
\hat{E}_{n}}$'s $\hat{E}_{n}$-surfaces.

Then we have the following results \cite{Leung-Xu-Zhang-12} generalizing the
construction of $E_{n}$-bundles over del Pezzo surfaces and their
relationships with flat $E_{n}$-bundles over elliptic curves.

\begin{Theorem}[\cite{Leung-Xu-Zhang-12}]
For any $\hat{E}_{n}$-surface $X_{\hat{E}_{n}}$, the orthogonal complement
of $K$ in $H^{2}\bigl( X_{\hat{E}_{n}},\mathbb{Z}\bigr) $ is isomorphic to
the root system of the affine Kac--Moody Lie algebra $\hat{E}_{n}$, i.e., \begin{equation*}
\langle K\rangle^{\bot }\simeq \Lambda_{\hat{E}_{n}}.
\end{equation*}%
Furthermore, the real root system of $\hat{E}_{n}$ is $\Delta^{\rm re}\bigl(
\hat{E}_{n}\bigr) \simeq \langle K\rangle^{\bot }\cap \{\alpha \cdot \alpha =-2\}$,
the imaginary root system of $\hat{E}_{n}$
is $\Delta^{\rm im}\bigl(\hat{E}_{n}\bigr)\simeq Z_{\neq 0}\langle K\rangle $ and the null root is $-K$.

There is also a canonically defined $\hat{E}_{n}$-bundle $\mathcal{E}^{\hat{E}_{n}}$ over $X_{\hat{E}_{n}}$,
\begin{align*}
\mathcal{E}^{\hat{E}_{n}} &=\mathcal{O}^{\oplus n+1}\oplus
\tbigoplus_{\alpha \in \Delta^{\rm re}(\hat{E}_{n})}O(\alpha)
\tbigoplus_{\beta \in \Delta^{\rm im}(\hat{E}_{n})}O(\beta)
\\
&\simeq \mathcal{E}^{E_{n}}\otimes \Bigl( \tbigoplus_{n\in \mathbb{Z}}O(nK) \Bigr) \oplus O.
\end{align*}
\end{Theorem}

Note that the last isomorphism is a realization of the Kac--Moody bundle as a
central extension of a loop algebra bundle over $X$, which depends on a
choice of a fixed line in $X$.

Given any fixed smooth elliptic curve $\Sigma $, let $\mathcal{S}_{\Sigma }^{%
\hat{E}_{n}}$ be the moduli space of $\hat{E}_{n}$-surfaces $X_{\hat{E}_{n}}$
containing $\Sigma $ as an anti-canonical curve and $\mathcal{M}_{\Sigma }^{%
\hat{E}_{n}}$ be the moduli space of holomorphic $\hat{E}_{n}$-bundles over $%
\Sigma $, in particular%
\begin{equation*}
\mathcal{M}_{\Sigma }^{\hat{E}_{n}}\simeq \operatorname{Hom}\bigl( \Lambda_{\hat{E}%
_{n}},\Sigma \bigr) /W_{\hat{E}_{n}\times \mathbb{Z}_{2}}.
\end{equation*}

\begin{Theorem}[\cite{Leung-Xu-Zhang-12}]
Given any smooth elliptic curve $\Sigma ,$ there is an open embedding
\begin{equation*}
\Phi\colon\ \mathcal{S}_{\Sigma }^{\hat{E}_{n}}\rightarrow \mathcal{M}_{\Sigma }^{\hat{E}_{n}}
\end{equation*}%
given by the restriction of the canonically defined $\hat{E}_{n}$-bundle $\mathcal{E}^{\hat{E}_{n}}$ over $X_{\hat{E}_{n}}$ to $\Sigma \subset X_{\hat{E}_{n}}$. Furthermore, there is a natural extension $\Phi\colon\mathcal{\bar{S}}_{\Sigma }^{\hat{E}_{n}}\overset{\sim }{\rightarrow }\mathcal{M}_{\Sigma }^{\hat{E}_{n}}$.
\end{Theorem}

In above discussions, we contract an $A_{d}$-chain of $(-2)$-curves on $X_{9}$ to realize affine $E_{9-d}$ structures on rational surfaces and elliptic curves. However, in Section~\ref{sec4}, $E_{n}$
structures on rational surfaces and elliptic curves also appear for
considerations of $\mathbb{P}^{2}$ blown up at $n$ points. In~general we
could consider $A_{d}$-chains of $(-2)$-curves on $\mathbb{P}^{2}$ blown up at $n+1$ points. Corresponding Lie algebras are listed in
the following magic triangle of Julia \cite{Julia}, see~Table~\ref{tab1}.

\begin{table}[!h]\centering\renewcommand{\arraystretch}{1.2}\setlength{\tabcolsep}{3.5pt}
\caption{Julia's magic triangle.}\label{tab1}\small \vspace{1mm}
\begin{tabular}{c|c|c|c|c|c|c|c|c|c}
\hline
& $n=8$ & $n=7$ & $n=6$ & $n=5$ & $n=4$ & $n=3$ & $n=2$ & $n=1$ & $n=0$
\\ \hline
$d=11$ & + &   &   &   &   &   &   &   &
\\ \hline
$d=10$ & $\mathbb{R}$ or $A_1$ & + &   &   &   &   &   &   &
\\ \hline
$d=9$ & $\mathbb{R}\times A_1$ & $\mathbb{R}$ &   &   &   &   &   &   &
\\ \hline
$d=8$ & $A_1\times A_2$ & $\mathbb{R}\times A_1$ or $A_2$ & $A_1$ &   &   &   &   &   &
\\ \hline
$d=7$ & $E_4$ & $\mathbb{R}\times A_2$ & $\mathbb{R}\times A_1$ & $\mathbb{R}$ & + &   &   &   &
\\ \hline
$d=6$ & $E_5$ & $A_1 \times A_3$ & $\mathbb{R}\times A_1^2$ & $\mathbb{R}^2$ or $A_1^2$ & $\mathbb{R}$ &   &   &   &
\\ \hline
$d=5$ & $E_6$ & $A_5$ & $A_2^2$ & $\mathbb{R}\times A_1^2$ & $\mathbb{R}\times A_1$ & $A_1$ &   &   &   \\ \hline
$d=4$ & $E_7$ & $D_6$ & $A_5$ & $A_1 \times A_3$ & $\mathbb{R}\times A_2$ & $\mathbb{R}\times A_1$ or $A_2$ & $\mathbb{R}$ & + &
\\ \hline
$d=3$ & $E_8$ & $E_7$ & $E_6$ & $E_5$ & $E_4$ & $A_1 \times A_2$ & $\mathbb{R}\times A_1$ & $\mathbb{R}$ or $A_1$ & +
\\  \hline
\end{tabular}
\end{table}

Notice that there is a symmetry between $d$ and $n$ in this magic triangle,
extending our earlier descriptions on $E_{n}$ bundles and affine $E_{n}$
bundles. The motivation of Julia is from the studies of 11-dimensional
supergravity in which the bosonic fields are pairs $(g,C)$,
where $g$ is a metric tensor on $\mathbb{R}^{1,10}$ and $C$ is a three form
field on $\mathbb{R}^{1,10}$, called the $C$-field. Consider a physical
toroidal compactification, namely one replaces $\mathbb{R}^{1,10}$ by $\mathbb{R}^{1,10-n}\times T^{n}$ and requires the volume of~$T^{n}$ to
shrink to zero size. By viewing $g$ as a family of metrics on $T^{n}$
parametrized by $\mathbb{R}^{1,10-n}$, we might expect to obtain an
effective theory which is a sigma model for maps from $\mathbb{R}^{1,10-n}$
to ${\rm SL}(n,\mathbb{Z}) \backslash {\rm SL}(n,\mathbb{R})
/{\rm SO}(n)$, as the latter space parametrizes Einstein metrics on $T^{n}$. However, the $C$-field will also decompose and regroup with components
of $g$ to enhance the sigma model to $E_{n,\mathbb{Z}}\backslash
E_{n}^{\rm split}/K_{n}$ where $K_{n}$ is the maximal compact subgroup of the
split Lie group $E_{n}^{\rm split}$ of type~$E_{n}$. Various structures that we
mentioned on ADE structures would have their counterparts in this
supergravity theory, which is a fascinating connection between algebraic
geometry and physics.

\section{ADE bundles over surfaces with ADE singularities}\label{sec8}
In this section, we explain that ADE singularities on a surface $X$
with $q=p_{g}=0$ leads to ADE bundles over $X$, as stated in Theorem~\ref{th8.1}.
Suppose $X^{\prime }$ is a singular surface with $q( X^{\prime})
=0$ and with a simple singularity $p$ and $X$ is its minimal resolution with
exceptional divisor $C=\bigcup_{i=1}^{r} C_{i}$. As explained in Section~\ref{sec2},
each irreducible component $C_{i}$ is an $(-2)$-curve, i.e., $C_{i}\simeq \mathbb{P}^{1}$ with $C_{i}\cdot C_{i}=-2$. In particular, $C_{i}\cdot K=0$ by the adjunction formula. The dual graph for the
configuration of these $C_{i}$'s is a Dynkin diagram of ADE type, thus there
is a corresponding simple Lie algebra $\mathfrak{g}$ of ADE type, and we
also call the corresponding singularity $p$ an ADE singularity. The $\mathbb{%
Z}$-span of the $C_{i}$'s is a root lattice $\Lambda_{\mathfrak{g}}$ of type $%
\mathfrak{g}$ inside $H^{2}(X,\mathbb{Z})$,
\begin{equation*}
\mathbb{Z}\langle C_{1},\dots,C_{r}\rangle =\Lambda_{\mathfrak{g}}\subset H^{2}(X,\mathbb{Z}).
\end{equation*}%
We denote the set of roots in $\Lambda_{\mathfrak{g}}$ as $\Phi $, i.e., $\alpha \in \Lambda_{\mathfrak{g}}$ lies in $\Phi $ if $\alpha \cdot \alpha=-2$. We write $\Phi =\Phi_{+}\amalg \Phi_{-}$, where $\alpha \in \Phi_{+}$ if $\alpha=\sum{n_iC_i}$ with $n_i\leq0$ and $\Phi_{-}=-\Phi_{+}$. Since $q(X) =0$, every class in $H^{2}(X,\mathbb{Z})\simeq \operatorname{Pic}(X)$ corresponds to a unique line bundle over $X$ up to isomorphisms.
Similar to our earlier constructions, we could construct a Lie algebra
bundle of ADE type $\mathfrak{g}$ over $X$ as follows:
\begin{equation*}
\mathcal{E}^{\mathfrak{g}}=O_{X}^{\oplus r}\tbigoplus\limits_{\alpha \in\Phi }O_{X}(\alpha).
\end{equation*}

It is natural to ask whether this ADE bundle $\mathcal{E}^{\mathfrak{g}}$
over $X$ can be descended to the original surface $X^{\prime }$ which admits
the ADE singularity $p$. Namely surfaces with an ADE singularity has a
natural ADE bundle over it. In \cite{Friedman-Morgan-02}, Friedman and
Morgan showed that it is possible for del Pezzo surfaces after small
deformations of $\mathcal{E}^{\mathfrak{g}}$. In \cite{Chen-Leung-14}, we
gave a direct construction of these deformations, which also works for
general surfaces $X^{\prime }$ satisfying $q=p_{g}=0$.

In the simplest case where $p$ is an $A_{1}$ singularity, namely locally $\mathbb{C}^{2}/\{\pm 1\}$, there is only one $(-2)$-curve $C_{1}$ in the exceptional locus and we have $\mathcal{E}^{A_{1}}=O_{X}\oplus O_{X}(-C_{1}) \oplus O_{X}(C_{1}) ={\rm End}_{0}(O_{X}\oplus O_{X}(C_{1}))$, the bundle of traceless endomorphisms of $O_{X}\oplus O_{X}(C_{1})$. Restricting to $C_{1}\subset X$ we have%
\begin{equation*}
(O_{X}\oplus O_{X}(C_{1})) |_{C_{1}}\simeq O_{\mathbb{P}^{1}}\oplus O_{\mathbb{P}^{1}}(-2),
\end{equation*}%
which admits a nontrivial deformation as an extension%
\begin{equation*}
0\rightarrow O_{\mathbb{P}^{1}}(-2) \rightarrow O_{\mathbb{P}^{1}}(-1) \oplus O_{\mathbb{P}^{1}}(-1) \rightarrow O_{\mathbb{P}^{1}}\rightarrow 0.
\end{equation*}%
Using $p_{g}(X) =0$, we could lift this extension from $C_{1}$
to $X$. The corresponding deformations of $\mathcal{E}^{A_{1}}$ would then
be trivial along $C_{1}$ and therefore can be descended to $X^{\prime }$ as
a Lie algebra bundle.

In the general case, we have the following result about these Lie algebra
bundles $\mathcal{E}^{\mathfrak{g}}$ through studying their minuscule
representation bundles in terms of $(-1)$-curves in $X$.

\begin{Theorem}[\cite{Chen-Leung-14}]\label{th8.1}
Let $p$ be an ADE singularity of a surface $X^{\prime }$ with $q(X^{\prime}) =p_{g}( X^{\prime }) =0$ and~$X$ be its minimal resolution with exceptional curve $C=\bigcup_{i=1}^{r} C_{i}$. Then
\begin{enumerate}\itemsep=0pt
\item[$(i)$] given any $(\varphi_{C_{i}})_{i=1}^{n}\in\Omega^{0,1}\big(X,\bigoplus
_{i=1}^{n}O(C_{i})\big)$ with $\overline{\partial}\varphi_{C_{i}}=0$ for every $i$, it can be extended to $\varphi=(\varphi_{\alpha})_{\alpha\in\Phi^{-}}\in\Omega^{0,1}\big(X,\bigoplus_{\alpha\in\Phi^{-}}O(\alpha)\big)$ such that $\overline{\partial}_{\varphi}:=\overline{\partial}+\operatorname{ad}(\varphi)$ is a holomorphic
structure on $\mathcal{E}^{\mathfrak{g}}$. We denote this new holomorphic
bundle as $\mathcal{E}_{\varphi}^{\mathfrak{g}}$;

\item[$(ii)$] such a $\overline{\partial}_{\varphi}$ is compatible with the Lie
algebra structure;

\item[$(iii)$] $\mathcal{E}_{\varphi}^{\mathfrak{g}}$ is trivial on $C_{i}$ if and
only if $[\varphi_{C_{i}}|_{C_{i}}]\neq0\in
H^{1}(C_{i},O_{C_{i}}(C_{i}))\cong\mathbb{C}$;

\item[$(iv)$] there exists $[\varphi_{C_{i}}]\in H^{1}(X,O(C_{i}))$ such that $[\varphi_{C_{i}}|_{C_{i}}]\neq0$;

\item[$(v)$] such a $\mathcal{E}_{\varphi}^{\mathfrak{g}}$ can descend to $X^{\prime }$ if and only if $[\varphi_{C_{i}}|_{C_{i}}]\neq0$ for every $i$.
\end{enumerate}
\end{Theorem}

\section{Relation to flag varieties of ADE type}\label{sec9}

In this section, we explain a relationship between our $G$-bundles
over a surface $X$ with a tautological $G$-bundle over the flag variety $G/B$
and its cotangent bundle  $G\times \mathfrak{n}/B$, as given in \cite[Theorems 5, 6 and 7]{Chen-Leung-18}.
An ADE singularity of type $\mathfrak{g}$ is locally given by the
intersection of the transversal slice $S_{x}$ of a subregular nilpotent
element $x$ and the nilpotent variety $N(\mathfrak{g})$ of the complex Lie
algebra $\mathfrak{g}$. Recall that $N(\mathfrak{g})$ is the fiber over zero
of the adjoint quotient $\mathfrak{g}\rightarrow \mathfrak{g/}G\simeq
\mathfrak{t}/W$, where $\mathfrak{t}$ is the Cartan subalgebra of $\mathfrak{g%
}$ and $W$ is the corresponding Weyl group. Furthermore, the restriction of
the adjoint quotient $\mathfrak{g}\rightarrow \mathfrak{t}/W$ to the
transversal slice~$S_{x}$ is a~semiuniversal deformation of the
corresponding $ADE$ singularity. This result is conjectured by Grothendieck
and proved by Brieskorn in 1970 \cite{Brieskorn}. After that, Grothendieck
defined a morphism $G\times \mathfrak{b}/B\rightarrow \mathfrak{t}$ and gave
a simultaneous resolution of the adjoint quotient $\mathfrak{g}\rightarrow
\mathfrak{t}/W$ using it. The restriction of the Grothendieck resolution to
the above transversal slice $S_{x}$ is also a simultaneous resolution \cite{Slodowy}. In 1969, Springer gave a resolution of singularities for the
nilpotent variety $N(\mathfrak{g})$ through $G\times \mathfrak{n}%
/B\rightarrow N(\mathfrak{g})$. Note that $G\times \mathfrak{n}/B\cong
T^{\ast }(G/B)$ is the cotangent bundle of the flag variety $G/B$. The
connection among these resolutions can be shown in the following
Brieskorn--Slodowy--Grothendieck diagram (here $\widetilde{S}$ is the minimal
resolution of $S$ and $C=\bigcup_{i=1}^{r} C_{i}$ is the exceptional locus
with each $C_{i}$ irreducible component):
\begin{equation*}
\begin{array}{l}
\begin{array}{*{20}c} {C=\bigcup_{i=1}^{r} C_{i}} & {\subset} &
{\widetilde{S}} & {\longrightarrow} & {S= N(\mathfrak{g}) \cap S_x } \\ \cap
& {} & \cap & {} & \cap \\ {G/B} & {\subset} & {G \times \mathfrak{n}/B} &
{\longrightarrow} & {N(\mathfrak{g})} \\ {} & {} & \cap & {\kern 1pt} & \cap
\\ {} & {} & {G \times \mathfrak{b}/B} & {\longrightarrow} & {\mathfrak{g}}
\\ {} & {} & \downarrow & {\kern 1pt} & \downarrow \\ {} & {} &
{\mathfrak{t}} & {\longrightarrow} & {\mathfrak{t}/W}. \\ \end{array} \\
\end{array}%
\end{equation*}

Under the above background, we consider the associated Lie algebra bundles $G\times \mathfrak{g}/B$ over $G/B$ and $G\times \mathfrak{n}\times \mathfrak{g}/B$ over $G\times \mathfrak{n}/B$ respectively. It is obvious that these
bundles are trivial as the action of $B$ on $\mathfrak{g}$ can extend to the
whole $G$. We \cite{Chen-Leung-18} describe natural holomorphic filtration
structures on these bundles explicitly.

As $B$ is a solvable Lie group, the associated representation bundle $%
G\times \mathfrak{g}/B$ over $G/B$ is an iterated extension of holomorphic
line bundles. Also for the full flag variety $G/B$, we have $%
\operatorname{Pic}(G/B)=\Lambda $, where $\Lambda $ is the weight lattice of the Lie
algebra $\mathfrak{g}$. Hence for every $\lambda \in \Lambda $, we can
associate a line bundle $L_{\lambda }$ over $G/B$. From the Borel--Weil--Bott
theorem, we can compute some particular cases of cohomology of line bundles
over $G/B$ easily. Note that our choice of $B$ gives rise to $\Phi_{+}$, the set of positive roots, such that
$\mathfrak{b} = \mathfrak{t} \oplus \bigoplus_{\alpha \in \Phi_{+}}\mathfrak{g}_{\alpha}$.
We let $\{\alpha_{1},\dots ,\alpha_{r}\}$ be the simple roots, then for any root $\alpha \in \Phi$, we have
\begin{enumerate}\itemsep=0pt
\item[$(i)$] $H^{i}(G/B,L_{\alpha })=0~~~$for any$~i\geq 2$;

\item[$(ii)$] $H^{1}(G/B,L_{\alpha }) \cong \mathbb{C}$ if $\alpha =-\alpha_{i}$ for some
simple root $\alpha_{i}$ and $H^{1}(G/B,L_{\alpha })=0$ otherwise;

\item[$(iii)$] the restriction map $H^{1}(G/B,L_{-\alpha_{i}})\rightarrow
H^{1}(C_{i},L_{-\alpha_{i}}|_{C_{i}})\cong \mathbb{C}$ is an isomorphism
for every simple root $\alpha_{i}$. Hence $[\varphi_{-\alpha
_{i}}|_{C_{i}}]\neq 0$ if and only if $[\varphi_{-\alpha_{i}}]\neq 0\in
H^{1}(G/B,L_{-\alpha_{i}})$.
\end{enumerate}

Now we try to write the holomorphic structures on $G\times \mathfrak{g}/B$
explicitly. The filtration of the representation $\mathfrak{g}$ is given by
the Chevalley order of its weights, hence not unique and we will choose an
arbitrary one. Then the holomorphic structure $\overline{\partial }_{\varphi
}$ on $G\times \mathfrak{g}/B$ can be written in an upper-triangular form
with respect to the holomorphic structure $\overline{\partial}$ on the graded
vector bundle associated to this filtration. Note that for a homogenous
space $G/P$, a vector bundle $V$ on $G/P$ is trivial if and only if the
restriction of $V$ to every Schubert line is trivial \cite{P}. Back to our
cases, the Schubert lines in $G/B$ are given by $C_{i}=P_{\alpha_{i}}/B$,
where $\alpha_{i}$'s run through all the simple roots, and $P_{\alpha_{i}}$
is the parabolic subgroup of $G$ corresponding to $\alpha_{i}$. The main
result is as follows:

\begin{Theorem}[\cite{Chen-Leung-18}]
For the Lie algebra bundle $(G\times \mathfrak{g}/B,\,[\, ,\,])$ over $G/B$ with
holomorphic structure $\overline{\partial }_{\varphi }$ as above, we have:
\begin{enumerate}\itemsep=0pt
\item[$(i)$] $\overline{\partial}_{\varphi} [\, ,\,]=0$ if and only if $\overline{\partial}_{\varphi}=\overline{\partial}+\sum_{\alpha\in\Phi^{-}}
    \operatorname{ad}(\varphi_{\alpha })$ with $\varphi_{\alpha} \in \Omega^{0,1}(G/B, L_{\alpha})$ for some $\alpha \in \Phi^-$.

\item[$(ii)$] The bundle $\bigl(G \times \mathfrak{g}/B, \overline{\partial}_{\varphi}=\overline{\partial}+\sum_{\alpha\in\Phi^{-}} \operatorname{ad}(\varphi_{\alpha })\bigr)$ is
holomorphically trivial if and only if $[\varphi_{-\alpha_i}|_{C_i}]\neq 0$
for every simple root $\alpha_i$.

\item[$(iii)$] The holomorphic structure of $(G \times \mathfrak{g}/B, \, [\, ,\,])$
over $G/B$ is $\overline{\partial}_{\varphi}=\overline{\partial}_{0} +\sum_{\alpha\in\Phi^{-}}\operatorname{ad}(\varphi_{\alpha })$ with $[\varphi_{-\alpha_i}]\neq 0$ for every simple root $\alpha_{i}$.
\end{enumerate}
\end{Theorem}

Consider the holomorphic structure of $G\times \mathfrak{n}\times \mathfrak{g}/B$ over $G\times \mathfrak{n}/B\cong T^{\ast }(G/B)$ when $\mathfrak{g}$
is of $ADE$ type. First, we know that $G\times \mathfrak{n}\times \mathfrak{g}/B$ is an iterated extension of line bundles over $G\times \mathfrak{n}/B$
as $B$ is solvable. And any line bundle over $G\times \mathfrak{n}/B$ is the
pull back of a line bundle over $G/B$ through the projection map $\pi
\colon T^{\ast }(G/B)\cong G\times \mathfrak{n}/B\rightarrow G/B$. Denote $\mathfrak{L}_{\lambda }:=\pi^{\ast }L_{\lambda }$ to be the corresponding
line bundle over $G\times \mathfrak{n}/B$ for any weight $\lambda \in
\Lambda $. Denote $H^{i}(\lambda ):=H^{i}(G\times \mathfrak{n}/B,\mathfrak{L}_{\lambda })$ for convenience. Using cohomology of line bundles on the
cotangent bundle of the flag variety \cite{Bro,Bro2,H}, we have:
\begin{enumerate}\itemsep=0pt
\item[$(i)$] for any positive root $\alpha \in \Phi^+$, $H^i(\alpha)=0$ for all $i\geq 1$;

\item[$(ii)$] for any negative root $\alpha \in \Phi^-$, $H^1(\alpha)\neq 0$, $H^2(\alpha)=0$;

\item[$(iii)$] the restriction map $H^{1}(G\times \mathfrak{n}/B, \mathfrak{L}_{-\alpha_i}) \rightarrow H^{1}(G/B, L_{-\alpha_i})$ is surjective for every
simple root $\alpha_{i}$.
\end{enumerate}

\begin{Theorem}[\cite{Chen-Leung-18}]
The holomorphic structure of $(G\times \mathfrak{n}\times \mathfrak{g}/B,\,[\, ,\,])$ over $G\times \mathfrak{n}/B$ is $\overline{\partial }_{\varphi
}=\overline{\partial }+\sum_{\alpha \in \Phi^{-}}\operatorname{ad}(\varphi_{\alpha })$
with $[\varphi_{-\alpha_{i}}|_{G/B}]\neq 0\in H^{1}(G/B,L_{-\alpha_{i}})$ for every simple root $\alpha_{i}$.
\end{Theorem}

Since the minimal resolution $\widetilde{S}$ of the $ADE$ singular surface $S$ is contained in $G\times \mathfrak{n}/B$, we also consider the
restriction of the $\mathfrak{g}$-bundle $G\times \mathfrak{n}\times
\mathfrak{g}/B$ from $G\times \mathfrak{n}/B$ to $\widetilde{S}$. Note that $\widetilde{S}$ is the minimal resolution of the $ADE$ singular surface $S=\mathbb{C}^{2}/\Gamma $. It is obvious that this $\mathfrak{g}$-bundle
over~$\widetilde{S}$ is also an iterated extension of line bundles. The
Picard group of $\widetilde{S}$ is a free abelian group generated by
divisors dual to the irreducible curves $C_{i}$ \cite{M}, i.e., $\operatorname{Pic}(\widetilde{S})=\mathbb{Z}\langle D_{i}\rangle $ with each $D_{i}$ dual to $C_{i}$.

As before, we know that the irreducible curves $C_{i}=P_{\alpha_{i}}/B$ are
Schubert lines in $G/B$, where~$\alpha_{i}$'s run through all the simple
roots. Now for any weight $\lambda $, we consider the restriction of the
line bundle $L_{\lambda }$ from $G/B$ to $C_{i}$ and it is easy to see that $L_{\lambda }|_{C_{i}}\cong \mathcal{O}_{\mathbb{P}^{1}}(\langle \lambda,\alpha_{i}\rangle )$. For the restriction of the line bundle $\mathfrak{L}_{\lambda }=\pi^{\ast }L_{\lambda }$ from $G\times \mathfrak{n}/B$ to $\widetilde{S}$, we know that for any root $\alpha =\sum n_{i}\alpha_{i}$, $\mathfrak{L}_{\alpha }|_{\widetilde{S}}\cong \mathcal{O}_{\widetilde{S}}\bigl(\sum -n_{i}C_{i}\bigr)$.

\begin{Theorem}[\cite{Chen-Leung-18}]
The restriction of the $\mathfrak{g}$-bundle $G\times \mathfrak{n}\times
\mathfrak{g}/B$ from $G\times \mathfrak{n}/B$ to $\widetilde{S}$ is
\begin{equation*}
\biggl(\mathcal{O}^{\oplus r}\oplus \bigoplus_{(\sum n_{i}C_{i})^{2}=-2}\mathcal{O}
\Bigl(\sum n_{i}C_{i}\Bigl),\,\overline{\partial }_{\varphi }=\overline{\partial}
+\sum_{\alpha \in \Phi^{-}}\operatorname{ad}(\varphi_{\alpha })\biggr)
\end{equation*}%
with $[\varphi_{-\alpha_{i}}]\neq 0\in H^{1}\big(\widetilde{S},\mathcal{O}%
(C_{i})\big)\cong \mathbb{C}$ for every simple root $\alpha_{i}$.
\end{Theorem}

We note that the holomorphic structures described here have the same form as
the holomorphic structures constructed in Section~\ref{sec8}.

\section{Generalization to affine ADE bundles}\label{sec10}

In this section, we will generalize the results for $ADE$ bundles in Section~\ref{sec8} to affine $ADE$ bundles and obtain \cite[Theorem~8]{Chen-Leung-16}. For convenience, we will call a curve $C=\cup
C_{i}$ in a surface $X$ an $ADE$ $($resp. affine $ADE)$ curve of type $%
\mathfrak{g}$ $($resp. $\widehat{\mathfrak{g}})$ if each $C_{i}$ is a smooth
$(-2)$-curve in $X$ and the dual graph of $C$ is a Dynkin diagram of the
corresponding type. By Kodaira's classification~\cite{Kodaira-63-I,Kodaira-63-II,Kodaira-63-III} of fibers of relative
minimal elliptic surfaces, every singular fiber is a~affine~$ADE$ curve
unless it is rational with a cusp, tacnode or triplepoint $($corresponding
to type~$II$ or $III\big(\widehat{A}_{1}\big)$ or $VI\big(\widehat{A}_{2}\big)$ in Kodaira's
notations$)$, which can also be regarded as a~degenerated affine $ADE$ curve
of type $\widehat{A}_{0}$, $\widehat{A}_{1}$ or $\widehat{A}_{2}$
respectively. We will not distinguish affine $ADE$ curves from their
degenerated forms since they have the same intersection matrices. We also
call the affine $ADE$ curves as Kodaira curves.

Suppose $C=\cup_{i=0}^{r}C_{i}$ is a affine $ADE$ curve of type $\widehat{%
\mathfrak{g}}$ in $X$ with $C_{0}$ corresponding to the extended root, then $%
\cup_{i=1}^{r}C_{i}$ will be the corresponding $ADE$ curve of type $\mathfrak{g}$ and
\begin{equation*}
\Phi :=\biggl\{\alpha =\Bigl[\sum_{i\neq 0}a_{i}C_{i}\Bigr]\in H^{2}(X,\mathbb{Z})\,|\,\alpha
^{2}=-2\biggr\}
\end{equation*}
is the root system of $\mathfrak{g}$. As before, we have a $\mathfrak{g}$-bundle
\begin{equation*}
\mathcal{E}^{\mathfrak{g}}=O^{\oplus r}\oplus \bigoplus_{\alpha \in \Phi}O(\alpha).
\end{equation*}
Also, there exists unique $n_{i}$'s up to overall scalings such that $%
F:=\sum n_{i}C_{i}$ satisfies $F\cdot F=0$. We know
\begin{equation*}
\Phi_{\widehat{\mathfrak{g}}}:=\{\alpha +nF\,|\,\alpha \in \Phi ,n\in \mathbb{Z}%
\}\cup \{nF\,|\,n\in \mathbb{Z},n\neq 0\}
\end{equation*}
is a affine root system and it decomposes into union of positive and
negative roots, i.e., $\Phi_{\widehat{\mathfrak{g}}}=\Phi_{\widehat{%
\mathfrak{g}}}^{+}\cup \Phi_{\widehat{\mathfrak{g}}}^{-}$, where
\begin{align*}
\Phi_{\widehat{\mathfrak{g}}}^{-}&=\Bigl\{\sum a_{i}C_{i}\in \Phi_{\widehat{\mathfrak{g}}}\,|\,a_{i}\geq 0 \text{ for all } i\Bigr\}
\\
&=\{\alpha +nF\,|\,\alpha \in \Phi^{+},\, n\in \mathbb{Z}_{\geq 1}\}\cup \{\alpha
+nF\,|\,\alpha \in \Phi^{-},\, n\in \mathbb{Z}_{\geq 0}\}\cup \{nF\,|\,n\in \mathbb{Z}_{\geq 1}\}
\end{align*}
and $\Phi_{\widehat{\mathfrak{g}}}^{+}=-\Phi_{\widehat{\mathfrak{g}}}^{-}$.
Then there is a canonically defined $\widehat{\mathfrak{g}}$-bundle $\mathcal{E}^{\widehat{\mathfrak{g}}}$ over $X$,
\begin{gather*}
\mathcal{E}^{\widehat{\mathfrak{g}}} =\mathcal{E}^{\mathfrak{g}}\otimes
\Bigl(\tbigoplus_{n\in \mathbb{Z}}O(nF)\Bigr) \oplus O.
\end{gather*}
We remark that if we remove the central extension part $O$ in $\mathcal{E}^{%
\widehat{\mathfrak{g}}}$, then this is a loop algebra $L\mathfrak{g}$ bundle
over $X$ and the Lie algebra bundle structure is independent of the choice
of the affine root $C_{0}$ in $F$. Similar to Section~\ref{sec8}, we have the
following results:

\begin{Theorem}[\cite{Chen-Leung-16}]
Given any complex surface $X$ with $p_{g}=0$. If $X$ has a Kodaira curve $%
C=\cup_{i=0}^{r}C_{i}$ of type $\widehat{\mathfrak{g}}$, then
\begin{enumerate}\itemsep=0pt
\item[$(i)$] given any $(\varphi_{C_{i}})_{i=0}^{r}\in \Omega^{0,1}\big(X,\bigoplus_{i=0}^{r}O(C_{i})\big)$ with $\overline{\partial}\varphi_{C_{i}}=0$ for every $i$, it can be extended to $\varphi =(\varphi_{\alpha})_{\alpha \in \Phi_{\widehat{\mathfrak{g}}}^{-}}\in \Omega^{0,1}\big(X,\bigoplus_{\alpha \in \Phi_{\widehat{\mathfrak{g}}}^{-}}O(\alpha)\big)$ such that $\overline{\partial}_{\varphi}:=\overline{\partial}+\operatorname{ad}(\varphi )$ is a holomorphic structure on $\mathcal{E}^{\widehat{\mathfrak{g}}}$. We denote the new bundle as $\mathcal{E}_{\varphi }^{\widehat{\mathfrak{g}}}$;

\item[$(ii)$] $\overline{\partial}_{\varphi}$ is compatible with the Lie algebra
structure on $\mathcal{E}^{\widehat{\mathfrak{g}}}$;

\item[$(iii)$] $\mathcal{E}_{\varphi}^{\widehat{\mathfrak{g}}}$ is trivial on $C_{i}$ if and only if $[\varphi_{C_{i}}|_{C_{i}}]\neq 0\in H^{1}(C_{i},O_{C_{i}}(C_{i}))\cong \mathbb{C}$;

\item[$(iv)$] there exists $[\varphi_{C_{i}}]\in H^{1}(X,O(C_{i}))$ such that $[\varphi_{C_{i}}|_{C_{i}}]\neq 0$.
\end{enumerate}
\end{Theorem}

\section[Deformability of Lie algebra bundles and geometry of rational surfaces]
{Deformability of Lie algebra bundles\\ and geometry of rational surfaces}\label{sec11}

In this section, we discuss the hidden geometry underlying the
deformability of the affine $E_{8}$ bundle over $X_{9}$, as stated in
\cite[Theorems 9 and~10]{Chen-Leung-15}. For an ADE curve $C$ in a surface $X$ with $p_{g}=0$, the corresponding Lie
algebra bundle $\mathcal{E}^{\mathfrak{g}}$ over $X$ admits a deformation
which can be descended to the surface obtained by blowing down $C$ in $X$.
On the other hand, an affine ADE curve can never be blown down. Nevertheless, we could explain the geometric meaning of such deformabilities as below.

Recall that we have an $E_{n}$ Lie algebra bundle $\mathcal{E}^{E_{n}}$ over
$X_{n}$ for $n\leq 8$. When $n=9$, $E_{9}$~is the affine Kac--Moody algebra
of $E_{8}$, i.e., $E_{9}=\hat{E}_{8}$. When $H^{2}(X,\mathbb{Z})$
has a sublattice $\Lambda_{\mathfrak{g}}$ isomorphic to the root lattice of
a simple Lie algebra $\mathfrak{g}$, then our construction also gives a $\mathfrak{g}$ Lie algebra bundle $\mathcal{E}^{\mathfrak{g}}$ over $X$.
Infinitesimal deformations of $\mathcal{E}^{\mathfrak{g}}$ as a $\mathfrak{g}$-bundle are parametrized by $H^{1}(X,\operatorname{ad}(\mathcal{E}^{\mathfrak{g}})) \simeq H^{1}(X,\mathcal{E}^{\mathfrak{g}})\subset H^{1}(X,\operatorname{End}(\mathcal{E}^{\mathfrak{g}}))$.
We say $\mathcal{E}$ is $(i)$ fully deformable if there is a~base $\Delta
\subset \Phi$ of the root system $\Lambda_g$ such that $H^{1}(X,O(\alpha)) \neq 0$ for every $\alpha \in \Delta$, $(ii)$ $\mathfrak{h}$-deformable with $\mathfrak{h}$ a Lie subalgebra of $\mathfrak{g}$ if there exists a strict $\mathfrak{h}$-subbundle of $\mathcal{E}$ which is fully deformable, $(iii)$
totally non-deformable if $H^{1}(X,O(\alpha)) =0$
for every $\alpha \in \Delta$, $(iv)$ deformable in $\alpha $-direction for $\alpha \in \Phi $ if $H^{1}(X,O(\alpha))\neq 0$.

We proved the following results.

\begin{Theorem}[\cite{Chen-Leung-15}]
On $X_{9}$, a blowup of $\mathbb{P}^{2}$ at $9$ points, these points are in
general position in $\mathbb{P}^{2}$ if and only if $\mathcal{E}^{E_{9}}$ is
totally non-deformable.
\end{Theorem}

\begin{Theorem}[\cite{Chen-Leung-15}]
If $-K$ is nef on $X_{9}$, then
\begin{enumerate}\itemsep=0pt
\item[$(i)$] there exists an ADE curve $C\subset X_{9}$ of type $\mathfrak{g}$ if and
only if $\mathcal{E}^{E_{9}}$ is $\mathfrak{g}$-deformable;

\item[$(ii)$] there exists a affine ADE curve $C\subset X_{9}$ of type $\widehat{\mathfrak{g}}$ if and only if $\mathcal{E}^{E_{9}}$ is $\widehat{\mathfrak{g}}$-deformable;

\item[$(iii)$] $X_{9}$ admits an elliptic fibration structure with a multiple fiber
of multiplicity $m$ if and only if $\mathcal{E}^{E_{9}}$ is deformable in the $(-mK)$-direction, but not in $(-m+1)K$-direction.
\end{enumerate}
\end{Theorem}

\section{Cox rings of ADE surfaces}\label{sec12}

In this section, we explain another mysterious relationship between del Pezzo surfaces $X_{n}$
and the Lie group $G$ of type $E_{n}$ discovered by Batyrev and Popov \cite{Batyrev-Popov-02} relating the Cox ring of $X_{n}$ and the flag variety $G/P_{L}\subset \mathbb{P}(L) $. We also establish the
the corresponding results for any ADE surface in~\cite[Theorem~11]{Leung-Zhang-15}. Very loosely speaking, the Cox
ring of a variety $X$ with $q(X) =0$ and torsion free $H^{2}(X,\mathbb{Z})$ is the sum of spaces of
sections of all line bundles on $X$,
\begin{equation*}
\operatorname{Cox}(X) \sim \tbigoplus\limits_{[L] \in \operatorname{Pic}(X) }H^{0}(X,L),
\end{equation*}%
with the ring structure given by the tensor products of sections. As
elements in $\operatorname{Pic}(X) $ are only isomorphism classes of line
bundles, in order to define $\operatorname{Cox}(X) $ properly, one needs to fix
a collection of line bundles $L_{i}$'s with $i=1,\dots ,b$ (with $%
b=b_{2}(X) $) whose first Chern classes represent a basis of $%
H^{2}(X,\mathbb{Z})$, then the correct definition of the Cox
ring of $X$ with respect to this choice is%
\begin{equation*}
\operatorname{Cox}(X) =\tbigoplus\limits_{n_{i}\in \mathbb{Z}}H^{0}\bigl(
X,L_{1}^{\otimes n_{1}}\otimes L_{2}^{\otimes n_{2}}\otimes \cdots \otimes
L_{b}^{\otimes n_{b}}\bigr) .
\end{equation*}%
Then the ring structure is simply given by the tensor products of
holomorphic sections.

When $X$ is a del Pezzo surface, the ample line bundle $K^{-1}$ defines a
grading on $\operatorname{Cox}(X) $. We have seen before that the geometry of a
degree $d$ del Pezzo surface is closely related to the Lie algebra $E_{n}$
with $n=9-d$. For example the number of lines in $X$ is equal to the dimension
of the fundamental representation $L$ for $n<8$. It turns out that the flag
variety $G/P_{L}$ corresponding to $L$, namely the unique closed $G$-orbit
in $\mathbb{P}(L) $, is closely related to the Cox ring of $X$
as follows: When the degree $d$ of $X$ satisfies $d\leq 5$, Batyrev--Popov
\cite{Batyrev-Popov-02}, Derenthal~\cite{Derenthal-07},
Serganova--Skorobogatov \cite{Serganova-Skorobogatov-07,Serganova-Skorobogatov-11} showed that there is a natural embedding of the projective
spectrum of the graded ring $\operatorname{Cox}(X) $ into $G/P_{L}$:
\begin{equation*}
\operatorname {Proj}(\operatorname{Cox}(X)) \hookrightarrow G/P_{L}.
\end{equation*}

It was observed in \cite{Leung-Zhang-15} that this relationship can be
easily generalized to all ADE surfaces~$X_{G}$. Recall that an ADE surfac~$X_{G}$ of rank $n$ is a blowup $X_{n+1}$ of $\mathbb{P}^{2}$ at $n+1$
distinct points, together with a rational curve $C$ in $X_{n+1}$ whose class
in $H^{2}(X_{n+1},\mathbb{Z})$ is $(i)$ $h$ satisfying $h\cdot
h=1 $ for type~$A$, $(ii)$ $f$ satisfying $f\cdot f=0$ for type $D$ and $(iii)$ $l$ satisfying $l\cdot l=-1$ for type~$E$. We~have $\langle
K,C\rangle^{\bot }\subset H^{2}(X_{n+1},\mathbb{Z})$ is
a root lattice $\Lambda_{G}$ of corresponding ADE type. In~particular, the
blow down of the $(-1)$-curve $l$ in an $E_{n}$-surface is a del
Pezzo surface $X$ of degree $d=9-n$ and $\operatorname{Cox}(X)
=\tbigoplus_{L\in \langle l\rangle^{\bot }}H^{0}(X_{En},L)$ defined loosely as before. Notice that $H^{2}(X,\mathbb{Z}) \simeq \langle l\rangle^{\bot }\subset H^{2}(X_{n+1},\mathbb{Z})$.

Similarly we define a generalization of the Cox ring for an ADE surface $%
X_{G}=(X_{n+1},C)$ of type $G$ as follows:
\begin{equation*}
\operatorname{Cox}^{G}(X_{G}) =\tbigoplus\limits_{L\in \langle C\rangle^{\bot }}H^{0}(X_{n+1},L).
\end{equation*}%
We also define the flag variety of $G$ corresponding to the fundamental
representation $L$ as $G/P_{L}$, i.e., $G/P_{L}$ is the unique closed $G$%
-orbit in $\mathbb{P}(L) $. Then we have the following statement:

\begin{Theorem}[\cite{Leung-Zhang-15}]
For any ADE surface $X_{G}=(X_{n+1},C)$ of type $G$, we have%
\begin{equation*}
\operatorname{Proj}\bigl(\operatorname{Cox}^{G}(X_{G})\bigr) \hookrightarrow G/P_{L},
\end{equation*}%
with $n\geq 4$ for type $E_{n}$ cases.
\end{Theorem}

\section{Conclusions}\label{sec13}

We have discussed several intriguing relationships between the geometry of
surfaces and Lie theory. A couple of these are related to physics, namely
the duality between F-theory and heterotic string theory and supergravity
theory in eleven dimensions. There is also a mysterious duality between the
geometry of del Pezzo surfaces and toroidal compactification of M-theory in
physics, as proposed by Iqbal, Neitzke and Vafa in \cite{Iqbal-Neitzke-Vafa}. We expect that more surprising connections will be uncovered in the future.

\subsection*{Acknowledgements}

Authors thank Jia Jin Zhang for many discussions and collaborations over
the years. And we also want to thank the referees for providing constructive
comments and help in improving the contents of this paper.
The work of N.~Leung described in this paper was substantially
supported by grants from the Research Grants Council of the Hong Kong
Special Administrative Region, China (Project No.~CUHK14301619 and
CUHK14306720) and a direct grant from CUHK (Project No.~CUHK4053400). The
work of Y.X.~Chen described in this paper was supported by Natural Science
Foundation of Shanghai (No.19ZR1411700).

\pdfbookmark[1]{References}{ref}
\LastPageEnding


\begin{thebibliography}{99}
\footnotesize\itemsep=0pt

\bibitem{BPV}
Barth W., Peters C., Van~de Ven A., Compact complex surfaces, \textit{Ergeb.
  Math. Grenzgeb.~(3)}, Vol.~4, \href{https://doi.org/10.1007/978-3-642-96754-2}{Springer-Verlag}, Berlin, 1984.

\bibitem{Batyrev-Popov-02}
Batyrev V.V., Popov O.N., The {C}ox ring of a del {P}ezzo surface, in
  Arithmetic of Higher-Dimensional Algebraic Varieties ({P}alo {A}lto, {CA},
  2002), \textit{Progr. Math.}, Vol.~226, \href{https://doi.org/10.1007/978-0-8176-8170-8_5}{Birkh\"auser Boston}, Boston, MA,
  2004, 85--103, \href{https://arxiv.org/abs/math.AG/0309111}{arXiv:math.AG/0309111}.

\bibitem{Bea}
Beauville A., Complex algebraic surfaces, 2nd ed., \textit{London Math. Soc. Stud.
  Texts}, Vol.~34, \href{https://doi.org/10.1017/CBO9780511623936}{Cambridge University Press}, Cambridge, 1996.

\bibitem{Beauville-99}
Beauville A., Counting rational curves on {$K3$} surfaces, \href{https://doi.org/10.1215/S0012-7094-99-09704-1}{\textit{Duke
  Math.~J.}} \textbf{97} (1999), 99--108, \href{https://arxiv.org/abs/alg-geom/9701019}{arXiv:alg-geom/9701019}.

\bibitem{Block-Gottsche-16}
Block F., G\"ottsche L., Refined curve counting with tropical geometry,
  \href{https://doi.org/10.1112/S0010437X1500754X}{\textit{Compos. Math.}} \textbf{152} (2016), 115--151, \href{https://arxiv.org/abs/1407.2901}{arXiv:1407.2901}.

\bibitem{Brieskorn}
Brieskorn E., Singular elements of semi-simple algebraic groups, in Actes du
  {C}ongr\`es {I}nternational des {M}ath\'ematiciens ({N}ice, 1970), {T}ome~2,
  Gauthier-Villars, Paris, 1971, 279--284.

\bibitem{Bro}
Broer B., Line bundles on the cotangent bundle of the flag variety,
  \href{https://doi.org/10.1007/BF01244299}{\textit{Invent. Math.}} \textbf{113} (1993), 1--20.

\bibitem{Bro2}
Broer B., Normality of some nilpotent varieties and cohomology of line bundles
  on the cotangent bundle of the flag variety, in Lie Theory and Geometry,
  \textit{Progr. Math.}, Vol.~123, \href{https://doi.org/10.1007/978-1-4612-0261-5_1}{Birkh\"auser Boston}, Boston, MA, 1994,
  1--19.

\bibitem{Bryan-Leung-Survey}
Bryan J., Leung N.C., Counting curves on irrational surfaces, in Surveys in
  Differential Geometry: Differential Geometry Inspired by String Theory,
  \textit{Surv. Differ. Geom.}, Vol.~5, \href{https://doi.org/10.4310/SDG.1999.v5.n1.a3}{Int. Press}, Boston, MA, 1999, 313--339.

\bibitem{Bryan-Leung-99}
Bryan J., Leung N.C., Generating functions for the number of curves on abelian
  surfaces, \href{https://doi.org/10.1215/S0012-7094-99-09911-8}{\textit{Duke Math.~J.}} \textbf{99} (1999), 311--328,
  \href{https://arxiv.org/abs/math.AG/9802125}{arXiv:math.AG/9802125}.

\bibitem{Bryan-Leung-YZ}
Bryan J., Leung N.C., The enumerative geometry of {$K3$} surfaces and modular
  forms, \href{https://doi.org/10.1090/S0894-0347-00-00326-X}{\textit{J.~Amer. Math. Soc.}} \textbf{13} (2000), 371--410,
  \href{https://arxiv.org/abs/alg-geom/9711031}{arXiv:alg-geom/9711031}.

\bibitem{Caporaso-Harris}
Caporaso L., Harris J., Counting plane curves of any genus, \href{https://doi.org/10.1007/s002220050208}{\textit{Invent.
  Math.}} \textbf{131} (1998), 345--392, \href{https://arxiv.org/abs/alg-geom/9608025}{arXiv:alg-geom/9608025}.

\bibitem{Chen-Leung-14}
Chen Y., Leung N.C., {$ADE$} bundles over surfaces with {$ADE$} singularities,
  \href{https://doi.org/10.1093/imrn/rnt065}{\textit{Int. Math. Res. Not.}} \textbf{2014} (2014), 4049--4084, \href{https://arxiv.org/abs/#2}{arXiv:1209.4979}.

\bibitem{Chen-Leung-15}
Chen Y., Leung N.C., Deformability of {L}ie algebra bundles and geometry of
  rational surfaces, \href{https://doi.org/10.1093/imrn/rnv023}{\textit{Int. Math. Res. Not.}} \textbf{2015} (2015),
  11508--11519.

\bibitem{Chen-Leung-16}
Chen Y., Leung N.C., Affine {$ADE$} bundles over surfaces with {$p_g=0$},
  \href{https://doi.org/10.1007/s00209-016-1645-6}{\textit{Math.~Z.}} \textbf{284} (2016), 55--68, \href{https://arxiv.org/abs/1303.5578}{arXiv:1303.5578}.

\bibitem{Chen-Leung-18}
Chen Y., Leung N.C., A{DE} bundles over {ADE} singular surfaces and flag
  varieties of {ADE} type, \href{https://doi.org/10.1093/imrn/rnx011}{\textit{Int. Math. Res. Not.}} \textbf{2018} (2018),
  3941--3958, \href{https://arxiv.org/abs/1811.02777}{arXiv:1811.02777}.

\bibitem{Clingher-Morgan-05}
Clingher A., Morgan J.W., Mathematics underlying the {F}-theory/heterotic
  string duality in eight dimensions, \href{https://doi.org/10.1007/s00220-004-1270-9}{\textit{Comm. Math. Phys.}} \textbf{254}
  (2005), 513--563, \href{https://arxiv.org/abs/math.AG/0308106}{arXiv:math.AG/0308106}.

\bibitem{Derenthal-07}
Derenthal U., Universal torsors of del {P}ezzo surfaces and homogeneous spaces,
  \href{https://doi.org/10.1016/j.aim.2007.01.012}{\textit{Adv. Math.}} \textbf{213} (2007), 849--864, \href{https://arxiv.org/abs/math.AG/0604195}{arXiv:math.AG/0604195}.

\bibitem{Dolgachev-2012}
Dolgachev I.V., Classical algebraic geometry: a modern view, \href{https://doi.org/10.1017/CBO9781139084437}{Cambridge
  University Press}, Cambridge, 2012.

\bibitem{Donagi-97}
Donagi R., Principal bundles on elliptic fibrations, \href{https://doi.org/10.4310/AJM.1997.v1.n2.a1}{\textit{Asian~J. Math.}}
  \textbf{1} (1997), 214--223, \href{https://arxiv.org/abs/alg-geom/9702002}{arXiv:alg-geom/9702002}.

\bibitem{Donagi-Wijnholt-20}
Donagi R., Wijnholt M., A{DE} transform, \href{https://doi.org/10.4310/ATMP.2020.v24.n8.a2}{\textit{Adv. Theor. Math. Phys.}}
  \textbf{24} (2020), 2043--2066, \href{https://arxiv.org/abs/1510.05025}{arXiv:1510.05025}.


\bibitem{M}
Donten-Bury M., Cox rings of minimal resolutions of surface quotient
  singularities, \href{https://doi.org/10.1017/S0017089515000221}{\textit{Glasg. Math.~J.}} \textbf{58} (2016), 325--355,
  \href{https://arxiv.org/abs/1301.2633}{arXiv:1301.2633}.

\bibitem{Fantechi-Gottsche-vanStraten-99}
Fantechi B., G\"ottsche L., van Straten D., Euler number of the compactified
  {J}acobian and multiplicity of rational curves, \textit{J.~Algebraic Geom.}
  \textbf{8} (1999), 115--133, \href{https://arxiv.org/abs/alg-geom/9708012}{arXiv:alg-geom/9708012}.

\bibitem{Friedman-Morgan-Witten}
Friedman R., Morgan J., Witten E., Vector bundles and {${\rm F}$} theory,
  \href{https://doi.org/10.1007/s002200050154}{\textit{Comm. Math. Phys.}} \textbf{187} (1997), 679--743,
  \href{https://arxiv.org/abs/hep-th/9701162}{arXiv:hep-th/9701162}.

\bibitem{Friedman-Morgan-02}
Friedman R., Morgan J.W., Exceptional groups and del {P}ezzo surfaces, in
  Symposium in {H}onor of {C}.{H}.~{C}lemens ({S}alt {L}ake {C}ity, {UT},
  2000), \textit{Contemp. Math.}, Vol.~312, \href{https://doi.org/10.1090/conm/312/04988}{Amer. Math. Soc.}, Providence, RI,
  2002, 101--116, \href{https://arxiv.org/abs/math.AG/0009155}{arXiv:math.AG/0009155}.

\bibitem{Gottsche-conj-98}
G\"ottsche L., A conjectural generating function for numbers of curves on
  surfaces, \href{https://doi.org/10.1007/s002200050434}{\textit{Comm. Math. Phys.}} \textbf{196} (1998), 523--533,
  \href{https://arxiv.org/abs/alg-geom/9711012}{arXiv:alg-geom/9711012}.

\bibitem{Gottsche-Shende-14}
G\"ottsche L., Shende V., Refined curve counting on complex surfaces,
  \href{https://doi.org/10.2140/gt.2014.18.2245}{\textit{Geom. Topol.}} \textbf{18} (2014), 2245--2307, \href{https://arxiv.org/abs/1208.1973}{arXiv:1208.1973}.

\bibitem{H}
Hague C., Cohomology of flag varieties and the {B}rylinski--{K}ostant
  filtration, \href{https://doi.org/10.1016/j.jalgebra.2009.04.001}{\textit{J.~Algebra}} \textbf{321} (2009), 3790--3815,
  \href{https://arxiv.org/abs/0803.3424}{arXiv:0803.3424}.

\bibitem{Hum}
Humphreys J.E., Introduction to {L}ie algebras and representation theory,
  \textit{Grad. Texts in Math.}, Vol.~9, \href{https://doi.org/10.1007/978-1-4612-6398-2}{Springer-Verlag}, New York~-- Berlin,
  1972.

\bibitem{Iqbal-Neitzke-Vafa}
Iqbal A., Neitzke A., Vafa C., A mysterious duality, \href{https://doi.org/10.4310/ATMP.2001.v5.n4.a5}{\textit{Adv. Theor. Math.
  Phys.}} \textbf{5} (2001), 769--807, \href{https://arxiv.org/abs/hep-th/0111068}{arXiv:hep-th/0111068}.

\bibitem{Julia}
Julia B.L., Magics of {${\rm M}$}-gravity, \href{https://doi.org/10.1002/1521-3978(200105)49:4/6<551::AID-PROP551>3.0.CO;2-Q}{\textit{Fortschr. Phys.}} \textbf{49}
  (2001), 551--555, \href{https://arxiv.org/abs/hep-th/0105031}{arXiv:hep-th/0105031}.

\bibitem{Kac}
Kac V.G., Infinite-dimensional {L}ie algebras, 3rd ed., \href{https://doi.org/10.1017/CBO9780511626234}{Cambridge University
  Press}, Cambridge, 1990.

\bibitem{KMPS-Pf-YZ-formula}
Klemm A., Maulik D., Pandharipande R., Scheidegger E., Noether--{L}efschetz
  theory and the {Y}au--{Z}aslow conjecture, \href{https://doi.org/10.1090/S0894-0347-2010-00672-8}{\textit{J. Amer. Math. Soc.}}
  \textbf{23} (2010), 1013--1040, \href{https://arxiv.org/abs/0807.2477}{arXiv:0807.2477}.

\bibitem{Kodaira-63-I}
Kodaira K., On compact analytic surfaces.~{I}, \href{https://doi.org/10.2307/1969881}{\textit{Ann. of Math.}}
  \textbf{71} (1960), 111--152.

\bibitem{Kodaira-63-II}
Kodaira K., On compact analytic surfaces.~{II}, \href{https://doi.org/10.2307/1970131}{\textit{Ann. of Math.}}
  \textbf{77} (1963), 563--626.

\bibitem{Kodaira-63-III}
Kodaira K., On compact analytic surfaces.~{III}, \href{https://doi.org/10.2307/1970500}{\textit{Ann. of Math.}}
  \textbf{78} (1963), 1--40.

\bibitem{Kool-Shende-Thomas-11}
Kool M., Shende V., Thomas R.P., A short proof of the {G}\"ottsche conjecture,
  \href{https://doi.org/10.2140/gt.2011.15.397}{\textit{Geom. Topol.}} \textbf{15} (2011), 397--406, \href{https://arxiv.org/abs/1010.3211}{arXiv:1010.3211}.

\bibitem{Lee-Leung-YZ-K3-nonprim}
Lee J., Leung N.C., Yau--{Z}aslow formula on {$K3$} surfaces for non-primitive
  classes, \href{https://doi.org/10.2140/gt.2005.9.1977}{\textit{Geom. Topol.}} \textbf{9} (2005), 1977--2012,
  \href{https://arxiv.org/abs/math.SG/0404537}{arXiv:math.SG/0404537}.

\bibitem{Lee-Leung-g=1}
Lee J., Leung N.C., Counting elliptic curves in {$K3$} surfaces,
  \href{https://doi.org/10.1090/S1056-3911-06-00439-5}{\textit{J.~Algebraic Geom.}} \textbf{15} (2006), 591--601,
  \href{https://arxiv.org/abs/math.SG/0405041}{arXiv:math.SG/0405041}.

\bibitem{Leung-ADE}
Leung N.C., {ADE}-bundles over rational surfaces, configuration of lines and
  rulings, \href{https://arxiv.org/abs/math.AG/0009192}{arXiv:math.AG/0009192}.

\bibitem{Leung-Xu-Zhang-12}
Leung N.C., Xu M., Zhang J., Kac--{M}oody {$\widetilde{E}_k$}-bundles over
  elliptic curves and del {P}ezzo surfaces with singularities of type~{$A$},
  \href{https://doi.org/10.1007/s00208-011-0661-4}{\textit{Math. Ann.}} \textbf{352} (2012), 805--828.

\bibitem{LEUNG-Zhang-ADE-I}
Leung N.C., Zhang J., Moduli of bundles over rational surfaces and elliptic
  curves. {I}.~{S}imply laced cases, \href{https://doi.org/10.1112/jlms/jdp053}{\textit{J.~Lond. Math. Soc.}} \textbf{80}
  (2009), 750--770, \href{https://arxiv.org/abs/0906.3900}{arXiv:0906.3900}.

\bibitem{LEUNG-Zhang-ADE-II}
Leung N.C., Zhang J., Moduli of bundles over rational surfaces and elliptic
  curves. {II}.~{N}onsimply laced cases, \href{https://doi.org/10.1093/imrn/rnp101}{\textit{Int. Math. Res. Not.}}
  \textbf{2009} (2009), 4597--4625, \href{https://arxiv.org/abs/0908.1645}{arXiv:0908.1645}.

\bibitem{Leung-Zhang-12}
Leung N.C., Zhang J., Non-simply laced {M}c{K}ay correspondence and triality,
  \href{https://doi.org/10.4310/PAMQ.2012.v8.n4.a6}{\textit{Pure Appl. Math.~Q.}} \textbf{8} (2012), 941--955.

\bibitem{Leung-Zhang-15}
Leung N.C., Zhang J., Cox rings of rational surfaces and flag varieties of
  {$ADE$}-types, \href{https://doi.org/10.4310/CAG.2015.v23.n2.a3}{\textit{Comm. Anal. Geom.}} \textbf{23} (2015), 293--317,
  \href{https://arxiv.org/abs/1409.2325}{arXiv:1409.2325}.

\bibitem{Manin-86}
Manin Yu.I., Cubic forms: algebra, geometry, arithmetic, 2nd ed., \textit{North-Holland
  Math. Libr.}, Vol.~4, North-Holland Publishing Co., Amsterdam, 1986.

\bibitem{P}
Pan X., Triviality and split of vector bundles on rationally connected
  varieties, \href{https://doi.org/10.4310/MRL.2015.v22.n2.a10}{\textit{Math. Res. Lett.}} \textbf{22} (2015), 529--547,
  \href{https://arxiv.org/abs/1310.5774}{arXiv:1310.5774}.

\bibitem{Pandharipande-Thomas-16}
Pandharipande R., Thomas R.P., The {K}atz--{K}lemm--{V}afa conjecture for {$K3$}
  surfaces, \href{https://doi.org/10.1017/fmp.2016.2}{\textit{Forum Math. Pi}} \textbf{4} (2016), e4, 111~pages,
  \href{https://arxiv.org/abs/1404.6698}{arXiv:1404.6698}.

\bibitem{Reid}
Reid M., Undergraduate algebraic geometry, \textit{London Math. Soc. Stud.
  Texts}, Vol.~12, \href{https://doi.org/10.1017/CBO9781139163699}{Cambridge University Press}, Cambridge, 1988.

\bibitem{Reid-2}
Reid M., La correspondance de {M}c{K}ay, \textit{Ast\'erisque} \textbf{276}
  (2002), 53--72, \href{https://arxiv.org/abs/math.AG/9911165}{arXiv:math.AG/9911165}.

\bibitem{Segre}
Segre B., The non-singular cubic surfaces, Oxford University Press, Oxford,
  1942.

\bibitem{Serganova-Skorobogatov-07}
Serganova V.V., Skorobogatov A.N., Del {P}ezzo surfaces and representation
  theory, \href{https://doi.org/10.2140/ant.2007.1.393}{\textit{Algebra Number Theory}} \textbf{1} (2007), 393--419,
  \href{https://arxiv.org/abs/math.AG/0611737}{arXiv:math.AG/0611737}.

\bibitem{Serganova-Skorobogatov-11}
Serganova V.V., Skorobogatov A.N., Adjoint representation of {${\rm E}_8$} and
  del {P}ezzo surfaces of degree~1, \href{https://doi.org/10.5802/aif.2676}{\textit{Ann. Inst. Fourier (Grenoble)}}
  \textbf{61} (2011), 2337--2360, \href{https://arxiv.org/abs/1005.1272}{arXiv:1005.1272}.

\bibitem{Slodowy}
Slodowy P., Simple singularities and simple algebraic groups, \textit{Lecture
  Notes in Math.}, Vol.~815, \href{https://doi.org/10.1007/BFb0090294}{Springer}, Berlin, 1980.

\bibitem{Tzeng-12}
Tzeng Y.-J., A proof of the {G}\"ottsche--{Y}au--{Z}aslow formula,
  \href{https://doi.org/10.4310/jdg/1335273391}{\textit{J.~Differential Geom.}} \textbf{90} (2012), 439--472,
  \href{https://arxiv.org/abs/1009.5371}{arXiv:1009.5371}.

\bibitem{Yau-Zaslow-96}
Yau S.-T., Zaslow E., B{PS} states, string duality, and nodal curves on {$K3$},
  \href{https://doi.org/10.1016/0550-3213(96)00176-9}{\textit{Nuclear Phys.~B}} \textbf{471} (1996), 503--512,
  \href{https://arxiv.org/abs/hep-th/9512121}{arXiv:hep-th/9512121}.

\end{thebibliography}
\end{document}